\documentclass{amsart}

\usepackage{graphicx}

\usepackage[utf8]{inputenc}
\usepackage{amssymb}

\newcommand{\field}{\mathbb{F}}

\newcommand{\prs}{\mathbb{P}}
\newcommand{\nat}{\mathbb{N}}

\newcommand{\nin}{\notin}
\newcommand{\union}{\cup}
\newcommand{\Union}{\bigcup}
\newcommand{\intersect}{\cap}
\newcommand{\Intersect}{\bigcap}
\newcommand{\id}{\operatorname{id}}

\newcommand{\defeq}{\mathrel{\mathop{:}}=}
\newcommand{\eqdef}{=\mathrel{\mathop{:}}}

\newcommand{\F}{\operatorname{F}}
\newcommand{\FP}{\operatorname{FP}}

\newcommand{\calA}{\mathcal{A}}
\newcommand{\calC}{\mathcal{C}}
\newcommand{\gen}[1]{\left\langle #1\right\rangle}
\newcommand{\tgen}[1]{{\langle {#1} \rangle}}
\newcommand{\abs}[1]{\left\lvert#1\right\rvert}
\newcommand{\chr}{\operatorname{char}}
\newcommand{\rad}{\operatorname{Rad}}
\newcommand{\typ}{\operatorname{typ}}
\newcommand{\pr}{\operatorname{pr}}
\newcommand{\proj}{\operatorname{proj}}
\newcommand{\coproj}{\operatorname{coproj}}
\newcommand{\realize}[1]{\abs{#1}}
\newcommand{\trealize}[1]{{\lvert {#1} \rvert}}
\newcommand{\st}{\operatorname{st}}
\newcommand{\lk}{\operatorname{lk}}
\newcommand{\SL}{\operatorname{SL}}
\newcommand{\SO}{\operatorname{SO}}
\newcommand{\Sp}{\operatorname{Sp}}
\newcommand{\Spin}{\operatorname{Spin}}
\newcommand{\stab}{\operatorname{Stab}}
\newcommand{\bfG}{\mathbf{G}}
\newcommand{\bfT}{\mathbf{T}}
\newcommand{\rsimeq}{\mathrel{\rotatebox[origin=c]{-90}{$\simeq$}}}

\newtheorem{thm}{Theorem}[section]
\newtheorem{lem}[thm]{Lemma}
\newtheorem{prop}[thm]{Proposition}
\newtheorem{cor}[thm]{Corollary}
\newtheorem*{maintheorem}{Main~Theorem}

\theoremstyle{definition}
\newtheorem{dfn}[thm]{Definition}

\theoremstyle{remark}
\newtheorem{rem}[thm]{Remark}

\numberwithin{equation}{section}

\newlength{\myimageoversize}
\newsavebox{\myimage}
\newcommand{\mycenter}[1]{%
\savebox{\myimage}{#1}
\settowidth{\myimageoversize}{\usebox{\myimage}}
\addtolength{\myimageoversize}{-\textwidth}
\setlength{\leftskip}{-0.5\myimageoversize}
\noindent
\usebox{\myimage}}

\begin{document}

\title[The sphericity of Phan geometries of type $B_n$ and $C_n$]{The sphericity of the Phan geometries \\ of type $B_n$ and $C_n$ \\ and \\ the Phan-type theorem of type $F_4$}

\author[Ralf Köhl]{Ralf Köhl (né Gramlich)}
\address{Mathematisches Institut, Universität Gießen, Arndtstraße 2, 35392 Gießen, Germany}
\email{ralf.koehl@math.uni-giessen.de}

\author{Stefan Witzel}
\address{Mathematisches Institut, Universität Münster, Einsteinstraße 62, 48149 Münster, Germany}
\email{s.witzel@uni-muenster.de}

\subjclass[2000]{Primary 51E24, Secondary 20G30, 20E42, 51A50}

\begin{abstract}
We adapt and refine the methods developed by Abramenko and De\-vil\-lers--Köhl--Mühlherr in order to establish the sphericity of the Phan geometries of type~$B_n$ and $C_n$, and their generalizations.

As an application we determine the finiteness length of the unitary form of certain hyperbolic Kac--Moody groups. We also reproduce the finiteness length of the unitary form of the groups $\Sp_{2n}(\field_{q^2}[t,t^{-1}])$.

Another application is the first published proof of the Phan-type theorem of type $F_4$. Within the revision of the classification of the finite simple groups this concludes the revision of Phan's theorems and their extension to the non-simply laced diagrams. We also reproduce the Phan-type theorems of types $B_n$ and $C_n$.
\end{abstract}

\maketitle

\section{Introduction}

In this paper we prove the following theorem; the exact statement is Theorem~\ref{thm:precise_statement} below, see Section~\ref{sec:gpgs} for definitions.

\begin{maintheorem}
A generalized Phan geometry of type $B_n$ or $C_n$ is $(n-1)$-spherical provided the defining field is sufficiently large. In fact, it is Cohen-Macaulay.
\end{maintheorem}

As the name suggests, the class of generalized Phan geometries contains the class of Phan geometries. These have been introduced in \cite{bengrahofshp07} (type~$B_n$) and \cite{grahofshp03} (type~$C_n$) in order to prove Phan-type theorems, that is, analogues of Phan's group-theoretic recognition results in \cite{phan77a,phan77b}. The Phan-type theorems state that the unitary forms of the groups $\Sp_{2n}(\field_{q^2})$ and $\mathrm{Spin}_{2n+1}(\field_{q^2})$ are the universal enveloping group of the amalgam of their fundamental subgroups of rank one and two. The unitary form is the subgroup fixed by the involution that acts as the field involution on $\field_{q^2}$ and takes $\field_q$-rational group elements to their transpose inverse. Using Tits' Lemma \cite[Corollaire 1]{tits86} this amalgamation result follows from the simple connectedness of the corresponding Phan geometries, which was established in \cite{bengrahofshp07} and \cite{grahornic07} in the case $B_n$ and in \cite{grahofshp03} and \cite{grahornic06} in the case $C_n$. This article provides alternative proofs of the Phan-type theorems of type $B_n$ and $C_n$ and, based on the local-to-global approach via the filtration described in \cite{devmue07}, gives a much shorter alternative to the unpublished proof of the Phan-type theorem of type $F_4$ by Hoffman, Köhl, M\"uhlherr, and Shpectorov. This concludes the revision of Phan's theorems and their extension to the non-simply laced diagrams, cf.\ \cite[page~656]{abrbro} and \cite[page~333]{gorlyosol05}.

Another class of geometries contained in the class of generalized Phan geometries are the geometries opposite one fixed chamber in a spherical building. The sphericity of these geometries (in spherical buildings of classical type) has been established by Abels and Abramenko in \cite{abeabr93}, \cite{abramenko96} in order to determine the topological finiteness length of certain arithmetic subgroups of affine Kac--Moody groups which are commensurable to the Borel lattice described in \cite{Carbone/Garland:1999}, \cite{Remy:1999}; cf.\ \cite{abels91}, \cite{abramenko96}.

Our strategy is an adaption of that work.
One of our main motivations for this article was to determine the finiteness length of the unitary form of certain affine Kac--Moody groups which are commensurable to the flip lattice described in \cite{gramue}. Concretely, if $\bfG$ is a semisimple algebraic group scheme of spherical type $X_n$, then $G \defeq \bfG(\field_{q^2}[t,t^{-1}])$ acts on a twin building $(\Delta_+,\Delta_-)$ of type $\tilde{X}_n$, see \cite[Chapter 11]{abrbro}. Let $\theta$ be the map that takes $g$ to $(g^\sigma)^{-T}$, where $\sigma$ acts as the field involution on $\field_{q^2}$ and exchanges $t$ and $t^{-1}$. We are interested in the group $K$ of $\theta$-fixed elements of $G$. The group $K$ is an arithmetic subgroup of $\bfG(\field_{q^2}(t))$ whose local rank over $\field_{q^2}((t))$ or, equivalently, $\field_{q^2}((t^{-1}))$ equals $n$. So, according to the Rank Conjecture, see \cite[p.~80]{behr98}, the group $K$ should be of finiteness type $F_{n-1}$, but not $FP_n$, which we confirm in the present article. Note that during the refereeing process for the present article the Rank Conjecture has been proved, see \cite{buxgrawit10b}.

The group $K$ naturally acts on the subcomplex $\Delta_\theta$ of the building $\Delta_+$ of chambers that are mapped to opposite ones by $\theta$. We use a filtration $(\Delta_i)_i$ of $\Delta_+$ with $\Delta_0 = \Delta_\theta$ from \cite{devmue07}, whose relative links will turn out (Theorem \ref{thm:flip_flops_are_gpgs}) to be generalized Phan geometries (cf.\ \cite[Fact 5.1]{devgramue08}). The types of generalized Phan geometries that occur are the types of spherical residues of the building. So by showing the sphericity of generalized Phan geometries, we show that the filtration preserves a certain degree of connectedness. Brown's criterion \cite[Corollary 3.3]{brown87} then allows us to determine the topological finiteness length of $K$.

Generalized Phan geometries of type $A_n$ have been investigated by Devillers, Köhl, and Mühlherr in \cite{devgramue08}. Since all irreducible spherical residues of a building of type $\tilde{A}_n$ are of type $A_m$, the procedure described above allowed the authors of loc.\ cit.\ to determine the finiteness length of the unitary form of $\SL_{n+1}(\field_{q^2}[t,t^{-1}])$. The Main Result of \cite{devgramue08} together with our Main Result allows us to determine the finiteness length of the unitary form of $\Sp_{2n}(\field_{q^2}[t,t^{-1}])$, see Theorem~\ref{thm:finiteness_length}. Both results provide an alternative proof of the Rank Theorem in their respective cases.

The method of proof also works for hyperbolic Kac--Moody groups provided all proper residues of the associated twin building are of spherical type, see Theorem~\ref{thm:finiteness_length_hyperbolic}.

The paper is organized as follows: The definition of generalized Phan geometries of type $A_n$, $B_n$ and $C_n$ is given in Section \ref{sec:gpgs}. In that section we also recall and collect some topological facts that are used later. Section \ref{sec:main_theorem} contains the precise statement of the Main Theorem and a proof of how it can be deduced from lemmas about the structure of certain subgeometries. These lemmas are then proved in Section \ref{sec:lemmas}. In Section \ref{sec:flip-flop-gpg} we will show that generalized Phan geometries occur as relative links in the filtration mentioned above. Finally, Section \ref{sec:applications} contains the applications described above, namely the computation of the topological finiteness length of the unitary form of $\Sp_{2n}(\field_{q^2}[t,t^{-1}])$ (Theorem~\ref{thm:finiteness_length}) and of certain hyperbolic Kac--Moody groups (Theorem~\ref{thm:finiteness_length_hyperbolic}), and the local recognition of groups that admit a weak Phan system of type $F_4$ (Theorem \ref{thm:local_recognition}).

\medskip
{\bf Acknowledgements:} The authors thank Bernhard M\"uhlherr for many very helpful discussions on how to properly define generalized Phan geometries. They also point out that the original unpublished proof of the Phan-type theorem of type $F_4$ was found by Hoffman, M\"uhlherr, Shpectorov and one of the authors during an RiP stay in Oberwolfach during the summer of 2005. Furthermore, we thank the referee for many helpful and valuable comments and suggestions.

\section{Basic definitions}
\label{sec:gpgs}

In this section we give the precise definition of the central objects of study of this paper, the generalized Phan geometries.

\medskip
Let $V$ be a vector space over a field $\field$. For $U,W \le V$, we say that $U$ is \emph{transversal} to $W$ and write $U \pitchfork W$, if $U \intersect W = 0$ or $\gen{U,W} = V$. Note that $U \pitchfork W$ if and only if $\dim(U \intersect W) = \max\{0,\dim U +\dim W - \dim V\}$.
For a flag $F=(0=V_0 \le \ldots \le V_k=V)$ and a subspace $U \le V$ we say that $U$ is \emph{transversal} to $F$ and write $U \pitchfork F$, if $U \pitchfork V_i$ for $0 \le i \le k$. A non-trivial $U$ is transversal to $F$ if and only if $\gen{U,V_{k_U}} = V$ where $k_U \defeq \min \{i \mid U \intersect V_j \ne \{0\}\}$.

Let $\sigma$ be an $\field$-automorphism of order one or two. Given a flag $F=(0=V_0 \le \ldots \le V_k=V)$ we call a family $(\omega_i)_{1 \le i \le k}$ of $\sigma$-hermitian forms $\omega_i \colon V_i \times V_i \to \field$ \emph{compatible with $F$} if $\rad \omega_i = V_{i-1}$.

\begin{rem}
Defining a $\sigma$-hermitian form on $V_i$ with radical $V_{i-1}$ is essentially the same as defining a non-degenerate $\sigma$-hermitian form on $V_i/V_{i-1}$.
\end{rem}

\begin{dfn}
Let $F$ and $\sigma$ be as above and let $\omega=(\omega_i)_i$ be a family of compatible $\sigma$-hermitian forms. For $U \le V$ we say that $U$ is \emph{transversal} to $(F,\omega)$, if $U$ is transversal to $F$ and $U \intersect V_{k_U}$ is $\omega_{k_U}$-non-degenerate. In this case we write $U \pitchfork (F,\omega)$.
\end{dfn}

We now define generalized Phan geometries of types $A_n$, $B_n$, and $C_n$. In the present work we are interested in types $B_n$ and $C_n$ but generalized Phan geometries of type $A_n$ will occur as join factors of residues, so we repeat the definition from \cite{devgramue08}:

\begin{dfn} \label{defan}
Let $F$, $V$, and $\sigma$ be as before and assume that $V$ has dimension $n+1$. Let
\[
F = (0 = V_0 \leq \ldots \leq V_k = V)
\]
be a flag an let $\omega = (\omega_i)_{1 \le i \le k}$ be a family of $\sigma$-hermitian forms that is compatible with $F$. The \emph{generalized Phan geometry of type $A_n$} defined by $(F,\omega)$ consists of all subspaces $U$ of $V$ that are transversal to $(F,\omega)$.
\end{dfn}

\begin{rem}
In contrast to \cite{devgramue08} we allow the inclusions of the $V_i$ to be non-strict. This is for purely technical reasons. If $V_{i-1} = V_i$ for some $i$, then one can clearly delete $V_i$ without changing the geometry. Note that in this case $\omega_i$ is the zero form.
\end{rem}

\begin{dfn} \label{defbn}
Assume that the characteristic of $\field$ is not $2$. Let $V$ be a vector space of dimension $2n+1$ equipped with a non-degenerate symmetric bilinear form $(\cdot,\cdot)$ of Witt index $n$. Let $F=(0=V_0 \le \ldots \le V_k=V)$ be a flag. Let $\omega$ be a family of $\sigma$-hermitian forms compatible with $F$ and assume that there is an $\omega_k$-non-isotropic vector that is $(\cdot,\cdot)$-isotropic. The \emph{generalized Phan-geometry of type $B_n$} defined by $(F,\omega)$ consists of all subspaces $U$ of $V$ that are totally $(\cdot,\cdot)$-isotropic and transversal to $(F,\omega)$.
\end{dfn}

\begin{dfn} \label{defcn}
Assume that the characteristic of $\field$ is not $2$ or that $\sigma \ne \id$. Let $V$ be a vector space of dimension $2n$ equipped with a non-degenerate alternating bilinear form $(\cdot,\cdot)$. Let $F=(0=V_0 \le \ldots \le V_k=V)$ be a flag. Let $\omega$ be a family of $\sigma$-hermitian forms compatible with $F$ and assume that there is an $\omega_k$-non-isotropic vector. The \emph{generalized Phan-geometry of type $C_n$} defined by $(F,\omega)$ consists of all subspaces $U$ of $V$ that are totally $(\cdot,\cdot)$-isotropic and transversal to $(F,\omega)$.
\end{dfn}

\begin{rem}
It is possible to define generalized Phan geometries using quadratic and unitary forms instead of bilinear and $\sigma$-hermitian forms, thus allowing to include characteristic $2$ in all cases. However, as $B_m(2^a) \cong C_m(2^a)$ by \cite[Theorem 2.2.10]{gorlyosol98} and since our main application concerns the case $\sigma \neq \mathrm{id}$, we refrain from adding even more technical difficulties to our study of generalized Phan geometries.
\end{rem}

\begin{dfn}
Let $V$ be a vector space equipped with a non-degenerate bilinear form. A family $(F^j)_j$ of flags of $V$ is \emph{matching} if for every $j$ there is an $l$ such that $(F^j)^\perp = F^l$. A family of generalized Phan geometries of same type $B_n$ or $C_n$ defined in a common vector space is matching if the flags defining them are matching.
\end{dfn}

Note that a single generalized Phan geometry defined by $(F,\omega)$ is matching if and only if $F = F^\perp$.

\medskip
We now recall and collect some topological facts, that we will need in the proof of the Main Theorem. For an introduction to simplicial complexes we refer the reader to \cite[Chapter 3]{spanier}.

In particular we use the following notions. If $K$ is a simplicial complex, $\trealize{K}$ is its \emph{realization}. If $s \in K$ is a simplex, the \emph{star of $s$} is $\st s = \{t \in K \mid t \union s \in K\}$. The \emph{link of $s$} is the subcomplex $\lk s = \{t \in \st s \mid t \intersect s = \emptyset\}$. The subcomplex generated by a subset $T \subseteq K$ is the complex $\overline{T} = \{s \in K \mid s \le t \textrm{ for some } t \in T\}$. We denote the boundary $\overline{\{s\}} \setminus \{s\}$ of a simplex $s$ by $\partial s$. If $L$ is another simplicial complex, then $K * L$ denotes the \emph{join} of the two, i.e. the complex with simplices $s \sqcup t, s \in K, t \in L$. This corresponds to the topological join in so far that $\trealize{K*L}$ is naturally homeomorphic to $\trealize{K}*\trealize{L}$.

An $n$-dimensional simplicial complex $K$ (as well as its realization $\trealize{K}$) is said to be \emph{$n$-spherical} if $\trealize{K}$ is $(n-1)$-connected. It is \emph{properly $n$-spherical}, if it is $n$-spherical not contractible. A simplicial complex is \emph{Cohen-Macaulay} if the link of every face is spherical. We will use the following facts, see \cite[Lemma~10.2, Theorem~10.6]{bjoerner95}:

\begin{prop}
\label{prop:topology_basics}
\begin{enumerate}
\item Let $X$ be a simplicial complex that can be written as $X = B \union \Union \calA$ where $B$ and $A \in \calA$ are subcomplexes. Assume that $\trealize{B}$ is $m$-connected, that $\trealize{A}$ is $m$-connected, that $\trealize{B \intersect A}$ is $(m-1)$-connected, and that $A \intersect A' \subseteq B$ for $A,A' \in \calA, A \ne A'$. Then $\trealize{X}$ is $m$-connected.\label{item:topological_filtration}
\item Let $X$ and $Y$ be simplicial complexes and assume that $\trealize{X}$ is $k$-connected and $\trealize{Y}$ is $l$-connected. Then $\trealize{X*Y}$ is $(k+l+1)$-connected.\label{item:join_of_sphericals}\qed
\end{enumerate}
\end{prop}

If $V$ is a vector space, there is a simplicial complex whose vertices are proper, nonempty subspaces of $V$ and whose simplices are flags of subspaces. If $\Gamma$ is any set of subspaces, then we may consider the full subcomplex on $\Gamma$, i.e.\ the subcomplex whose vertex set is $\Gamma$ and whose simplices are flags of elements of $\Gamma$. We will sometimes say that $\Gamma$ has a certain property to mean that the associated simplicial complex has that property.

We also need to fix some notation concerning buildings and twin buildings. Since we study topological properties of buildings and their sub-geometries, we consider buildings as simplicial complexes as introduced for example in \cite[Chapter 4]{abrbro}. 

If $\Delta$ is a building of type $I$, then the chamber systems of $\st s$ and $\lk s$ are isomorphic and describe the residue of type $I \setminus \typ(s)$ of any chamber in $\st s$. However $\lk s$ is the more natural model in so far as it is the simplicial complex naturally obtained from the chamber system. On the other hand $\st s$ has the virtue that its chambers actually are (and not only correspond to) chambers of $K$. This simplifies notation when dealing with projections. We will therefore consider both complexes as residues and make the correspondence implicit.

\section{The Main Theorem}
\label{sec:main_theorem}

In this section we show how the following theorem, which is the exact statement of our Main Theorem, can be deduced from four lemmas concerning the structure of certain subgeometries. The lemmas will be proven in the following section.

\begin{thm}
\label{thm:precise_statement}
Let $m,n > 0$ be integers, let $\field$ be a field, and let $\sigma$ be an automorphism of $\field$ of order $1$ or $2$. If $\field$ is finite, assume that $\abs{\field} \ge 4^{n-1}2 m$, if $\sigma = \id,$ and that $\abs{\field} \ge 4^{n-1}(q+1)m$, if $\field=\field_{q^2}$ and $\sigma \neq \id$.

Let $V$ be an $\field$-vector space and let $(\Gamma^j)_{1 \le j \le m}$ be a matching family of generalized Phan geometries in $V$ with respect to $\sigma$, all of same type either $B_n$ or $C_n$. Then the intersection $\Gamma = \Intersect^j \Gamma_j$ is $(n-1)$-spherical.
\end{thm}

Let $V$ be a vector space equipped with a bilinear form $(\cdot,\cdot)$ and let $\Gamma$ be the geometry of the theorem given by a flags $F^j=(0=V_0^j < \ldots < V_{k}^j=V)$ and families $\omega^j=(\omega_i^j)_i$ of $\sigma$-hermitian forms $\omega_i^j \colon V_i^j \times V_i^j \to \field$.

We use the following strategy due to Abramenko \cite[Chapter II]{abramenko96} to prove the theorem. Let $p \in \Gamma$ be one-dimensional.
Let $Z=\{U \in \Gamma \mid \gen{p,U} \pitchfork (F^j,\omega^j), 1 \le j \le m\}$. Let $Y_0 = \{U \in Z \mid U \intersect p^\perp,\gen{U\intersect p^\perp,p} \in \Gamma\}$.
For $1 \le i \le n$ let
\[
Y_i=Y_{i-1} \union \{U \in Z \mid \dim(U)=i\}
\]
and for $n+1 \le i \le 2n$ let
\[
Y_i=Y_{i-1} \union \{U \in \Gamma \mid \dim(U)=2n+1-i\}\textrm{ .}
\]

We are going to show that each $Y_i$ is $(n-1)$-spherical by induction on $i$. For this purpose it suffices to show that the residue of $U \in Y_i \setminus Y_{i-1}$ in $Y_i$ is $(n-2)$-spherical: we can then apply Proposition~\ref{prop:topology_basics}~\eqref{item:topological_filtration} with $X = Y_i$, $B = Y_{i-1}$, and $\calA = \{\st_{Y_i} U \mid U \in Y_i \setminus Y_{i-1}\}$. Sphericity of the residues will follow if we can prove the lemmas below.

To understand the general idea let us look at $U \in \Gamma$ of dimension $k$. The full residue of $U$ in $\Gamma$ consists of all spaces $W \in \Gamma$ such that either $W < U$ or $W > U$. Hence it is the join
\[
\{W < U \mid W \in \Gamma\} * \{W > U \mid W \in \Gamma\}\text{ .}
\]
It is not hard to see that the left join factor is an intersection of generalized Phan geometries of type $A_{k-1}$ and the right part is an intersection of generalized Phan geometries of type $C_{n-k}$ (the latter statement is Lemma~\ref{lem:upper_residue_is_cn} below). To understand which connectedness properties are preserved by the passage from $Y_{i-1}$ to $Y_i$ one has to study the subgeometries of the residue of $U \in Y_i$ of spaces $W$ with $W \in Y_{i-1}$. The filtration is constructed in such a way that these are again intersections of generalized Phan geometries. This is what the lemmas state.

\begin{lem}
\label{lem:upper_residue_intersected_is_cn}
If $U \in Y_i \setminus Y_{i-1}$ for $1 \le i \le n$, then $Y_i^{>U} \defeq \{W \in Y_i \mid W > U\} = \{W \in Y_0 \mid W > U\}$ is an intersection of at most $4m$ matching generalized Phan geometries of type $B_{n - k}$ respectively $C_{n - k}$.
\end{lem}

\begin{lem}
\label{lem:upper_residue_is_cn}
If $U \in Y_i \setminus Y_{i-1}$ for $n+1 \le i \le 2n$, then $Y_i^{>U}=\{W \in Y_i \mid W > U\} = \{W \in \Gamma \mid W > U\}$ is an intersection of at most $m$ matching generalized Phan geometries of type $B_{n-k}$ respectively $C_{n - k}$.
\end{lem}

\begin{lem}
\label{lem:lower_residue_is_an}
The set $Y_i^{<U} \defeq \{W \in Y_i \mid W < U\} = \{W \in Z \mid W < U\}$ is an intersection of at most $2m$ generalized Phan geometries of type $A_{k-1}$.
\end{lem}

\begin{lem}
\label{lem:induction_basis}
Let $\Gamma$ be the intersection of $m$ generalized Phan geometries, all of type either $B_1$ or $C_1$, defined over $\field$ by $\sigma$-hermitian forms. If $\field$ is finite, assume that $\abs{\field} \ge 2m$, if $\sigma = \id$, and that $\abs{\field} \ge (q+1)m$, if $\sigma \neq \id$ and $\field=\field_{q^2}$. Then $\Gamma$ contains at least one point.
\end{lem}

\begin{proof}[Proof of Theorem \ref{thm:precise_statement}]
To simplify notation define $C := 2$ in case $\sigma = \id$ and $C := q+1$ in case $\sigma \neq \id$, so that the condition on the field in the statement of the theorem can be written as $\abs{\field} \ge 4^{n-1}Cm$.

We proceed by an induction on $n$. To be $0$-spherical, the space has to be $(-1)$-connected, i.e. non-empty, so the induction basis is just Lemma \ref{lem:induction_basis}.
Assume that the statement is true for all $k < n$ and consider the filtration $(Y_i)_i$ described above.

We now proceed by an induction on $i$.
The geometry $Y_0$ is contractible because $U \mapsto U \intersect p^\perp \mapsto \gen{U \intersect p^\perp, p} \mapsto p$ is a deformation retraction to one point.

We want to apply Proposition~\ref{prop:topology_basics}~\eqref{item:topological_filtration} to the setup $X = Y_i$, $B = Y_{i-1}$, and $\calA = \{\st_{Y_i}(U) \mid U \in Y_i \setminus Y_{i-1}\}$. The space $\trealize{Y_{i-1}}$ is $(n-1)$-spherical by induction. To see that $\trealize{\st_{Y_i}(U)} \intersect \trealize{\st_{Y_i}(U')} \subseteq \trealize{Y_{i-1}}$, let $W \in \st_{Y_i}(U) \intersect \st_{Y_i}(U')$, i.e., $W$ is incident with $U$ and $U'$. Then its dimension is not $\dim(U) = \dim(U')$, whence $W \in Y_{i-1}$. The space $\trealize{\st_{Y_i}(U)}$ is clearly contractible, as it is a cone over $\trealize{\lk_{Y_i}(U)}$. So it remains to see that, for $U \in Y_i \setminus Y_{i-1}$ of dimension $k$, the space $ \trealize{\st(U) \intersect Y_{i-1}}$ is $(n-2)$-spherical.

In order to do so, we remark that $\st(U) \intersect Y_{i-1} = Y_i^{<U} * Y_i^{>U}$. By Lemma \ref{lem:lower_residue_is_an}, the complex $Y_i^{<U}$ is the intersection of at most $2m$ generalized Phan geometries of type $A_{k-1}$. Since $2^{k-1}C2m \le 4^{n-1}Cm \le q$, the Main Theorem from \cite{devgramue08} implies that $Y_i^{<U}$ is $(k-1)$-spherical. Similarly, Lemmas \ref{lem:upper_residue_intersected_is_cn} or \ref{lem:upper_residue_is_cn} imply that $Y_i^{>U}$ is an intersection of at most $4m$ matching generalized Phan geometries of type $B_{n-k}$ or $C_{n-k}$. Now $4^{k-2}C4m \le 4^{n-1}Cm \le q$, so by induction $Y_i^{>U}$ is $(n-k-1)$-spherical. Hence $Y_i^{<U} * Y_i^{>U}$ is $(n-1)$-spherical by Proposition~\ref{prop:topology_basics}~\eqref{item:join_of_sphericals}.
\end{proof}

\section{Proof of the lemmas}
\label{sec:lemmas}

Before proceeding to the proof of the lemmas let us first recall some notions from linear algebra. Let $V$ be an $\field$-vector space, $\sigma$ be an $\field$-automorphism of order $1$ or $2$, and $\varepsilon \in \{-1,1\}$. A form $(\cdot,\cdot) \colon V \times V \to \field$ is called {\em $(\sigma,\varepsilon)$-hermitian}, if it is biadditive and $(\lambda x, \mu y)=\lambda \mu^\sigma (x,y)$ and $(x,y)=\varepsilon (y,x)^\sigma$ for all $x,y \in V$, $\lambda,\mu \in \field$, $x,y \in V$. It is {\em $\sigma$-hermitian}, if it is $(\sigma,1)$-hermitian. Clearly symmetric bilinear, $\sigma$-hermitian and alternating forms are $(\sigma,\varepsilon)$-hermitian for some $\sigma$, $\varepsilon$. 
The form is \emph{non-degenerate} if $V^\perp=0$.

\begin{lem}
\label{lem:la}
Let $V$ be a vector space and let $A$, $B$ and $C$ be subspaces.
\begin{enumerate}
\item Then $\dim(\gen{A,B})=\dim(A)+\dim(B)-\dim(A \intersect B)$.\label{item:la_dim}
\item Let $(\cdot,\cdot)$ be a $(\sigma,\varepsilon)$-hermitian form. Then $A^\perp \intersect B^\perp = \gen{A,B}^\perp$. If in addition $(\cdot,\cdot)$ is non-degenerate, then also $\dim A^\perp = \dim V - \dim A$, $A^{\perp\perp}=A$,
and $\gen{A^\perp,B^\perp} = (A \intersect B)^\perp$.\label{item:la_forms}\qed
\end{enumerate}
\end{lem}

\begin{lem}
\label{lem:perp_in_quotient}
Let $V$ be a vector space that is equipped with a $(\sigma,\varepsilon)$-hermitian form. If $U \le V$ is isotropic, then $U^\perp/U$ naturally carries a quotient $(\sigma,\varepsilon)$-hermitian form. This form satisfies $W^\perp/U = (W/U)^{\perp'}$ for every $U \le W \le U^\perp$.\qed
\end{lem}

%

We start with the proof of the induction basis, Lemma \ref{lem:induction_basis}. Here is a restatement (in fact a slightly stronger version) that is obtained by unfolding the definitions.

\begin{lem}[restatement of Lemma \ref{lem:induction_basis}]
Let $\field$ be a field and
\begin{enumerate}
\item let $V$ be an $\field$-vector space of dimension $2$ and $(\cdot,\cdot)$ be a non-degenerate alternating bilinear form on $V$ (case $C_1$), or
\item let $V$ be an $\field$-vector space of dimension $3$ and $(\cdot,\cdot)$ be a non-degenerate symmetric bilinear form of Witt index $1$ on $V$ (case $B_1$).
\end{enumerate}
Let $\omega^j$, $1 \le j \le m$, be $\sigma$-hermitian forms on $V$. Assume that $\chr \field \ne 2$ unless we are in case $C_1$ and $\sigma \ne \id$. If $\field$ is finite, assume that $\abs{\field} > 2m$, if $\sigma = \id$, and $q^2 > (q+1)m$, if $\field=\field_{q^2}$ and $\sigma \neq \id$.
Assume that for every $j$ there is a vector that is $(\cdot,\cdot)$-isotropic and $\omega^j$-non-isotropic. Then there is a vector that is $(\cdot,\cdot)$-isotropic and $\omega^j$-non-isotropic for $1 \le j \le m$ simultaneously.
\end{lem}

\begin{proof}
We proceed by induction on the number of forms. The induction basis holds by hypothesis. Assume the statement to be true up to $m-1$.

\emph{Case $C_1$:}
First we consider the alternating case. The condition that a vector be $(\cdot,\cdot)$-isotropic is empty, so we are left with the condition on the forms $\omega^j$. Assume that $\sigma=\id$ and let $y,z$ form a basis for $V$.
If $y$ is non-isotropic with respect to all of the $\omega^j$, then we are done. Otherwise, consider the vectors $u_\alpha = \alpha y + z$ with $\alpha \in \field$. Then
\begin{equation}
\omega^j(u_\alpha,u_\alpha) = \alpha^2 \omega^j(y,y) + \omega^j(z,z) + 2\alpha \omega^j(y,z) \textrm{ ,} \quad 1 \le j \le m,
\label{eq:quadratic_form}
\end{equation}
are non-zero polynomials in $\alpha$ that are at most quadratic, and at least one of the polynomials (one with $\omega^j(y,y) = 0$) is at most linear. Each polynomial has at most two zeroes (and those with $\omega^j(y,y) = 0$ have at most one zero), so if $\abs{\field} \ge 2m$, there exists an $\alpha$ such that $\omega^j(u_\alpha,u_\alpha) \ne 0$ for all $i$.

Now assume that $\sigma \neq \id$. If $\field$ is infinite, let $\field^\sigma$ denote the field of its $\sigma$-fixed elements. Let $V'$ be the $\field^\sigma$-span of an $\field$-basis of $V$. Then $V'$ is a $2$-dimensional $\field^\sigma$-vector space on which the $\omega^j$ are symmetric bilinear. By the argument above there exists a vector $u \in V'$ that is $\omega^j$-non-isotropic.

If $\field = \field_{q^2}$ and $\sigma \neq \id$, then $\alpha^\sigma=\alpha^q$ for $\alpha \in \field$. The argument now works as above with \eqref{eq:quadratic_form} replaced by
\begin{equation}
\omega^j(u_\alpha,u_\alpha) = \alpha^{q+1}\omega^j(y,y) + \omega^j(z,z) + (\alpha + \alpha^q)\omega^j(y,z)\textrm{ ,} \quad 1 \le j \le m \textrm{ .}
\end{equation}
The condition becomes that $\abs{\field} \ge (q+1)m$.

\emph{Case $B_1$:}
Now assume that $(\cdot,\cdot)$ is a non-degenerate symmetric bilinear form of Witt index $1$ on a vector space of dimension $3$. By \cite[Theorem 6.3.1]{cameron91} the space $(V,(\cdot,\cdot))$ is the direct sum of a $2$-dimensional hyperbolic space and a $1$-dimensional non-degenerate space. The same is true of $(V,\xi(\cdot,\cdot))$ where $\xi$ is a non-square. Since we are only interested in whether or not vectors are $(\cdot,\cdot)$-isotropic, it does not matter which of the two forms we consider. One of them -- we assume without loss of generality, that it be $(\cdot,\cdot)$ -- has Gram matrix
\[
\left(
\begin{array}{ccc}
&1&\\
1&&\\
&&1
\end{array}
\right)
\]
with respect to a basis $e,f,x$. The vectors $f$ and $e - \beta^2/2 f + \beta x$ with $\beta \in \field$ are up to scalar multiples all isotropic vectors of $(\cdot,\cdot)$.

Again we consider first the case where $\sigma = \id$. If $e$ is non-isotropic with respect to all $\omega^j$, then we are done. Otherwise, we consider $u = e- \beta^2/2 f + \beta x$ for $\beta \in \field$. As before, $\omega^j(u,u)$ is a non-zero polynomial of degree at most $2$ in $\beta$, so it has at most $2$ zeroes. And if $\omega^j(e,e) = 0$ then it has at most one zero. As before we find that if $\abs{\field} \ge 2m$ we find a vector that is $(\cdot,\cdot)$-isotropic but $\omega^j$-non-isotropic for $1 \le i \le m$.

The case $\sigma \neq \id$ can be derived from the $\sigma = \id$-case just as in the $C_1$-case.
\end{proof}

Next we prove Lemmas \ref{lem:upper_residue_intersected_is_cn}, \ref{lem:upper_residue_is_cn}, and \ref{lem:lower_residue_is_an}.
The proofs consists essentially in considering the case of an individual generalized Phan geometry (i.e.\ $m=1$) and constructing one, two or four generalized Phan geometries that define the join factors of the residues. For arbitrary $m$ in the cases $B_n$ and $C_n$ it remains then to verify, that the $m$ respectively $4m$ generalized Phan geometries obtained in this way are matching if the original ones were. This is important because of

\begin{lem}
\label{lem:perp_is_transversal}
Let $(F^j)_{j}$ be a matching family of flags and let $U$ be transversal to all $F^j$. Then $U^\perp$ is also transversal to all $F^j$.
\end{lem}

\begin{proof}
Let $j$ and $l$ be such that ${V_i^j}^\perp = V_{t-i}^l$. By Lemma~\ref{lem:la}~\eqref{item:la_forms} $U + V_i^j = V$ is equivalent to $U^\perp \intersect V_{t-i}^l = 0$ and $U \intersect V_i^j = 0$ is equivalent to $U^\perp + V_{t-i}^l = V$.
\end{proof}

Lemma \ref{lem:lower_residue_is_an} is essentially contained in \cite{devgramue08}:

\begin{lem}
\label{lem:flag_for_gen_p_intersect}
Let $F = (0=V_0 < \ldots < V_k = V)$ be a flag and let $\omega=(\omega_i)_i$ be $\sigma$-hermitian forms compatible with $F$. Let $U \pitchfork F$ and let $p$ be a one-dimensional $\omega_k$-non-degenerate space such that $p \not \le U$. Then there is a flag $F'$ of $U$ and a family of forms $\omega'=(\omega_i')_i$ compatible with $F'$ such that
\[
\gen{p,W} \pitchfork (F,\omega) \quad \textrm{ if and only if } \quad W \pitchfork (F',\omega')
\]
for $W < U$.
\end{lem}

\begin{proof}
This can be read off the proof of \cite[Lemma 4.11]{devgramue08}.
\end{proof}

\begin{proof}[Proof of Lemma \ref{lem:lower_residue_is_an}]
Since $U$ is totally isotropic with respect to $(\cdot,\cdot)$, any $W < U$ is also totally isotropic with respect to $(\cdot,\cdot)$. Hence $\{W \in Z \mid W < U\}$ consists of all $W < U$ such that $W \pitchfork (F^j,\omega^j)$ and $\gen{p,W} \pitchfork (F^j,\omega^j)$ for all $j$. Let $F_1^j=(V_i \intersect U)_{m \le i \le k}$, $\omega_1^j=(\omega_i|_{U})_{m \le i \le k}$ where $m = \max\{i \mid V_i \intersect U=0\}$. Moreover, let $(F_2^j,\omega_2^j)$ be the restriction to $U$ of the flag whose existence is guaranteed by Lemma \ref{lem:flag_for_gen_p_intersect}. Then $\{W \in Z \mid W < U\}$ is the intersection of the $2m$ generalized Phan geometries defined on $U$ by $(F_1^j,\omega_1)$ and $(F_2^j,\omega_2)$, $1 \le j \le m$.
\end{proof}

We proceed to the proof of Lemma \ref{lem:upper_residue_is_cn}.

That is, we consider an isotropic subspace $U \le V$ of dimension $k$ that is transversal to $(F,\omega)$ and have to show that the geometry of isotropic subspaces $W > U$ that are transversal to $(F,\omega)$ form a geometry that is isomorphic to a generalized Phan geometry of type $B_{n-k}$ or $C_{n-k}$.

Note that every such $W$ is contained in $U' \defeq U^\perp$. The vector space in which we construct the generalized Phan geometry is $U'/U$ and the isomorphism will be $W \mapsto W/U$. Note that $(\cdot,\cdot)$ naturally induces a form on $U'/U$ (by Lemma~\ref{lem:perp_in_quotient} applied to $U'$) and that $W/U$ is totally isotropic if and only if $W$ is.

Let a generalized Phan geometry on $V$ be given by a flag $F = (0 = V_0 < \ldots < V_k = V)$ and forms $\omega = (\omega_i)_i$. The data by which the generalized Phan geometry on $U'/U$ is defined is constructed as follows:

Let $m = \max \{i \mid \gen{V_i,U} \intersect U' = U\}$ and let $M = \min \{i \mid \gen{V_i,U} \intersect U' = U'\}$. Let $A = (U \intersect V_M)^{\perp_{\omega_M}} \intersect U'$. We define the flag $F'$ in $U'/U$ by $F' = (V'_m \le \cdots \le V'_M)$ with $V_i' = \gen{V_i \intersect U',U}/U$. The forms $\omega' = (\omega_i')_{m \le i \le M}$ are given by $\omega_i'(x+U,y+U) = \omega_i(x,y)$ for $x,y \in A$.

Recall that we write $k_W \defeq \min\{i \mid V_i \intersect W \ne 0\}$ for any subspace $W$ of $V$. Similarly we write $k'_{W/U} \defeq \min\{i \mid V_i' \intersect W/U \ne 0\} = \min\{i \mid V_i \intersect W \not \le U\}$ for any subspace $W/U$ of $U/U$. The same symbols will be used throughout the section; the flag to which they refer should be clear from the context.

\begin{lem}
\label{lem:flag_for_project}
Let $V$ be a vector space, let $F$ be a flag in $V$ and let $\omega$ be a family of $\sigma$-hermitian forms compatible with $F$. Let $U$ be a subspace of $V$ that is transversal to $(F,\omega)$ and let $U' \ge U$ be transversal to $F$. The flag $F'$ and the forms $\omega'$ described above have the properties that $\omega'$ is compatible with $F'$ and if $W$ is such that $U \le W \le U'$, then
\[
W \pitchfork (F,\omega) \quad \textrm{ if and only if }\quad W/U \pitchfork (F',\omega') \text{ .}
\]
\end{lem}

\begin{proof}
We may assume $U' \gneq U$.

Note that $M = \min\{i \mid \tgen{V_i, U} = V\}$: If $\tgen{V_i,U} = V$, then trivially $\gen{V_i,U} \intersect U' = U'$. Conversely, if $\gen{V_i,U} \intersect U' = U'$, then $\tgen{V_i,U} \ge U',V_i$. But $U' \pitchfork F$, so $V_i \intersect U' \ne 0$ (which holds because $U' \gneq U$) implies $\gen{V_i,U} \ge \tgen{U',V_i} = V$. As a consequence $k_U=M$ or $k_U=M+1$. Similarly, using $U \pitchfork F$ one sees that $m = \max\{i \mid V_i \intersect U' = 0\}$, hence $m = k_{U'} - 1$.

Note further that $A$ is a complement for $U \intersect V_{M}$ in $U'$ because $U \intersect V_{M}$ is $\omega_{M}$-non-degenerate and $\tgen{U,V_{M}} \intersect U' = U'$. Note also that $U' \intersect V_i = (A \intersect V_i) \oplus_{\omega_i} (U \intersect V_i)$ for every $m \le i \le M$: indeed, either $U \intersect V_i = 0$  and $V_i \intersect U' \le \rad \omega_{k_U} \intersect U' \le A$ or $i = k_U$ and $U' \intersect V_{k_U} = (U \intersect V_{k_U}) \oplus_{\omega_{k_U}} A$.

We show first that $W$ is transversal to $F$ if and only if $W/U$ is transversal to $F'$. If $i \ge M$, i.e., $\gen{U,V_i} = V$, then $\gen{W,V_i} = V$ and $\gen{W,(V_i+U) \intersect U'} = U'$. Similarly, if $i \le m$, i.e. $U' \intersect V_i = 0$, then $W \intersect V_i = 0$ and $W \intersect \gen{V_i,U} = U$. So it suffices, indeed, to restrict attention to the cases $m \le i < M$.

If $W \intersect V_i = 0$, then $W/U \intersect \gen{V_i, U}/U = U/U$. If $\gen{W,V_i} = V$, then also $\gen{W/U,(\gen{V_i,U} \intersect U')/U} = U'/U$. Conversely, if $W/U \intersect \gen{V_i,U}/U = U/U$, then $W \intersect V_i \le U \intersect V_i = 0$. So assume that $W/U \intersect (\gen{V_i,U} \intersect U')/U \ne U/U$ and $\gen{W/U, (\gen{V_i,U} \intersect U')/U} = U'/U$. Then $\gen{V_i,W} \ge V_i,U'$ and, as $U'$ is transversal to $V_i$ and $U' \intersect V_i \ne 0$, we get that $\gen{U',V_i} = V$. The preceding discussion shows also that $k_W = k'_{W/U}$.

Now we show that the forms $\omega_i'$ are compatible with $F'$. To show that $\rad \omega_i' \le V_{i-1}'$ let $a \in A \intersect V_i$ be such that $a + U \in \rad \omega_i'$ and let $y \in V_i$ be arbitrary. We can write $y = b + u$ with $b \in A \intersect V_i$, $u \in U \intersect V_i$ and have
\[
0 = \omega_i'(a+U,b+U) = \omega_i(a,b) = \omega_i(a,b) + \omega_i(a,u) = \omega_i(a,y) \textrm{ ,}
\]
so $a \in \rad \omega_i = V_{i-1}$ and, thus, $a + U \in V_{i-1}'$.

Conversely, assume that $x + U \in V_{i-1}'$ and write $x = a + u$ for $a \in A \intersect V_i$, $u \in U \intersect V_i$. Then, for $b \in A \intersect V_i$, we have
\[
0 = \omega_i(x,b) = \omega_i(a,b) + \omega_i(u,b) = \omega_i(a,b) = \omega_i'(a+U,b+U) \textrm{ ,}
\]
whence $x + U = a + U \in \rad \omega_i'$.

Finally we show that $W \intersect V_{k_W}$ is $\omega_{k_W}$-non-degenerate if and only if $W/U \intersect V_{k_W}'$ is $\omega_{k_W}'$-non-degenerate.
Let $x \in (W \intersect V_{k_W}) \intersect (W \intersect V_{k_W})^{\perp_{\omega_{k_W}}}$. We write $x = u + a$ with $u \in U \intersect V_{k_W}$ and $a \in A \intersect V_{k_W}$. For $y \in W \intersect V_{k_W} \intersect A$ we then have
\[
\omega_{k_W}'(x + U, y + U) = \omega_{k_W}(a,y) = \omega_{k_W}(x,y) = 0
\]
where the second equality holds because $u \in U \intersect V_{k_W}$ and $y \in A \intersect V_{k_W}$. Hence $x + U \in \rad \omega_{k_W}'$.

Conversely, assume that $x \in W \intersect V_{k_W} \intersect A$ is such that $\omega_{k_W}'(x+U,y+U) = 0$ for all $y \in W \intersect V_{k_W} \intersect A$. Then, for $z \in W \intersect V_{k_W}$ and writing $z = u + w$ with $u \in U \intersect V_{k_W}$ and $y \in V_{k_W} \intersect A$, we get
\[
\omega_{k_W}(x,z) = \omega_{k_W}(x,u) + \omega_{k_W}(x,y) = \omega_{k_W}'(x+U,y+U) = 0 \text{ ,}
\]
since $x \in V_{k_W} \intersect A$ and $u \in V_{k_W} \intersect U$. Hence $x \in \rad \omega_{k_W}$.
\end{proof}

\begin{proof}[Proof of Lemma \ref{lem:upper_residue_is_cn}]
Let $(\Gamma^j)_{1 \le j \le m}$ be a family of matching generalized Phan geometries of type $B_n$ or $C_n$ defined by a flags $(F^j)$ and families of forms $(\omega^j)_{j}$. Let $\Gamma \defeq \Intersect_{1 \le j \le m} \Gamma^j$ and let $U \in \Gamma$. For every $j$ let ${\Gamma^j}'$ be the generalized Phan geometry on $U'/U$ with flag ${F^j}'$ and forms ${\omega^j}'$ obtained from $\Gamma^j$ as in Lemma \ref{lem:flag_for_project}. Then the geometry of $\{ W < V \mid U < W < U'\}$ is isomorphic to $\Gamma' \defeq \Intersect_{1 \le j \le m} {\Gamma^j}'$. Moreover, if $\Gamma^l$ is such that $(F^j)^\perp = F^l$, then ${{F^j}'}^\perp = {F^l}'$ by Lemma~\ref{lem:perp_in_quotient}.
\end{proof}

For the proof of Lemma \ref{lem:upper_residue_intersected_is_cn} we will need a somewhat refined version of Lemma \ref{lem:flag_for_project}. Namely, we have to weaken the assumptions that $U$ and $U'$ be transversal to $F$ and instead allow deviations by one dimension.

\begin{dfn}
A space $U$ is said to be \emph{almost transversal} to $F$, if $U \intersect V_i \ne 0$ implies that $\gen{U,V_i}$ has codimension at most $1$. And $U$ is said to be \emph{almost transversal} to $(F,\omega)$, if it is almost transversal and $U \intersect V_{k_U}$ is $\omega_{k_U}$-non-degenerate.
A space $U$ is said to be \emph{nearly transversal} to $F$, if $\gen{U,V_i} \ne V$ implies that $U \intersect V_i$ has dimension at most $1$.
\end{dfn}

Consider a flag and forms $(F,\omega)$ and subspace $U \lneq U' \le V$ such that $U$ is almost transversal to $(F,\omega)$ and $U'$ is nearly transversal to $F$. We want to show that the geometry of subspaces $W$ with $U < W < U'$ that are transversal to $(F,\omega)$ form a geometry that is isomorphic to a geometry of subspaces of $U'/U$ transversal to an appropriately chosen $(F',\omega')$.

\begin{lem}
\label{lem:generalized_projection}
Let $V$ be a vector space, $F = (0=V_0 \le \ldots \le V_k=V)$ be a flag and $\omega=(\omega_i)_i$ be a family of $\sigma$-hermitian forms compatible with $F$. Let $U \le U' \le V$ and assume that $U$ is almost transversal to $(F,\omega)$ and $U'$ is nearly transversal to $F$. Let $m = \max \{i \mid \gen{V_i,U} \intersect U' = U\}$ and $M = \min\{i \mid \gen{V_i,U} \intersect U' = U'\}$. There is an $A$ such that $U' \intersect V_i = A \intersect V_i \oplus_{\omega_i} U \intersect V_i$.

Consider the flag $F' = (V_{m}' \le \ldots \le V_{M}')$ with $V_i' = (\gen{V_i,U} \intersect U')/U$. The forms $\omega'=(\omega_i')_i$ given by $\omega_i'(x+U,y+U) = \omega_i(x,y)$ for $x,y \in A$ are compatible with $F'$ and 
\[
W \pitchfork (F,\omega) \quad \textrm{ if and only if }\quad W/U \pitchfork (F',\omega')
\]
for $U < W < U'$.
\end{lem}

\begin{proof}
Note that, if $\dim U' - \dim U \le 1$, the statement is trivial because there is no $W$ with $U < W < U'$. So we may and do assume for the rest of the proof that $\dim U' - \dim U \ge 2$.

The points at which we used that $U \pitchfork F$ in the proof of Lemma~\ref{lem:flag_for_project} are the following:
\begin{enumerate}
\item the construction of $A$,
\item $m = \max \{i \mid V_i \intersect U' = 0\}$, \label{item:m_u_characterization}
\item if $W/U \intersect (\gen{V_i,U} \intersect U')/U = U/U$, then $W \pitchfork F$\label{item:W_quotient_trivial}.
\end{enumerate}

The following construction of $A$ works without any assumptions on the transversality of $U$ to $F$. Let $A_1$ be an $\omega_1$-orthogonal complement for $U \intersect V_1$ in $U' \intersect V_1$. Note that $A_1 \le \rad \omega_2$, so $A_1$ is $\omega_2$-orthogonal to $U \intersect V_2$ and can be extended to an $\omega_2$-orthogonal complement of $U \intersect V_2$ in $U' \intersect V_2$. Iterating this process we get finally that $A = A_{M}$ is as desired.

To see that \eqref{item:m_u_characterization} is still true under the new hypotheses, assume that $\gen{V_i,U} \intersect U' = U$. Then $V_i \intersect U' \le U$. If $V_i \intersect U \ne 0$, then $\gen{V_i,U}$ would have codimension at most $1$, hence $\gen{V_i,U} \intersect U'$ would have codimension at most $1$ in $U'$. But $U$ has codimension at least $2$ in $U'$, hence $\gen{V_i,U} \intersect U' \ne U$ contradicting the assumption. So we see that in fact $V_i \intersect U = 0$.

For \eqref{item:W_quotient_trivial} we have to show that, if $W \intersect \gen{V_i,U} = U$, then $W$ is transversal to $F$. As before we have that $W \intersect V_i = U \intersect V_i$. If $W \intersect V_i = 0$, then there is nothing to show, so we may assume $U \intersect V_i \ne 0$ (i.e. $i = k_U$). Then, since $U$ is almost transversal to $F$, we know that $\gen{U,V_i}$ has codimension at most $1$. So $W \gneq U$ and $W \intersect \gen{V_i,U} = U$ imply $\gen{V_i,W} = V$.


The points at which we used that $U' \pitchfork F$ in the proof of Lemma~\ref{lem:flag_for_project} are the following:
\begin{enumerate}
\item $M = \min \{i \mid \gen{V_i, U} = V\}$, \label{item:M_u_characterization}
\item if $W/U \intersect (\gen{V_i,U}/U) \ne U/U$ and $\gen{W/U,\gen{V_i,U}/U} = U'/U$, then $\gen{W,V_i} = V$\label{item:W_span_all}.
\end{enumerate}

To see that \eqref{item:M_u_characterization} is still true under the new hypotheses, assume that $U' = \gen{V_i,U} \intersect U' =  \gen{U,V_i \intersect U'}$. Since $\dim U' - \dim U > 1$, we see that $V_i \intersect U'$ must have dimension at least two. So $\gen{V_i,U'} = V$, because $U'$ is nearly transversal to $F$.

For \eqref{item:W_span_all} assume that $W/U \intersect \gen{V_i,U}/U \ne U$ and $\gen{W/U,\gen{V_i,U}/U} = U'/U$. Then $\gen{W,V_i} \ge V_i,U'$, i.e.\ $\gen{W,V_i}=\gen{U',V_i}$. Hence $U' \intersect V_i \gneq W \intersect V_i \ne 0$. This shows that $\dim (U' \intersect V_i) \ge 2$ so that $\gen{U',V_i} = V$ because $U'$ is nearly transversal to $F$.
\end{proof}

The following is elementary:

\begin{lem}
\label{lem:isomorphisms}
Let $V$ be a finite-dimensional vector space, and let $U < U'$ be subspaces. Let $H$ be a hyperplane of $V$ such that $U \not \le H$. There is an isomorphism of vector spaces $\phi \colon U'/U \to (U' \intersect H)/(U \intersect H)$ such that
\[
\phi(W/U) = (W \intersect H)/(U \intersect H)
\]
for $U \le W \le U'$.
Similarly, let $p$ be a one-dimensional subspace of $V$ such that $p \not \le U'$. 

There is an isomorphism of vector spaces $\phi \colon U'/U \to \gen{U',p}/\gen{U,p}$ such that
\[
\phi(W/U) = \gen{W,p}/\gen{U,p}
\]
for $U \le W \le U'$.\qed
\end{lem}

\begin{proof}[Proof of Lemma \ref{lem:upper_residue_intersected_is_cn}]
Note that $U \le p^\perp$ and $U \in Z$ imply $U \in Y_0$. Hence $U \not \le p^\perp$. For the same reason, every $W$ in
\begin{align*}
Y_0^{>U} & = \{W \in Y_0 \mid W > U\}\\
& = \{U < W < U^\perp \mid W \in \Gamma, W \intersect p^\perp \in \Gamma, \gen{W \intersect p^\perp,p} \in \Gamma, \gen{p,W} \pitchfork (F,\omega)\}
\end{align*}
satisfies $W \not \le p^\perp$.
Note that $W \le U^\perp$ implies that $W \intersect p^\perp$ is totally isotropic if and only if $W = \gen{W \intersect p^\perp,U}$ is totally isotropic. Furthermore, note that $W \intersect p^\perp$ is totally isotropic if and only if $\gen{W \intersect p^\perp,p}$ is totally isotropic. Hence $Y_0^{>U}$ consists of the spaces $W$ with $U< W < U^\perp $ that are totally isotropic and satisfy
\begin{enumerate}
\item $W \pitchfork (F^j,\omega^j)$,
\item $\gen{W,p} \pitchfork (F^j,\omega^j)$,
\item $W \intersect p^\perp \pitchfork (F^j,\omega^j)$,
\item $\gen{W \intersect p^\perp,p} \pitchfork (F^j,\omega^j)$
\end{enumerate}
for $1 \le j \le m$.

Let us collect some properties of the spaces mentioned above. Clearly, $U$ and $\gen{U,p}$ are transversal to $(F^j,\omega^j), 1 \le j \le m$, because $U \in Z$. So by Lemma \ref{lem:perp_is_transversal} the spaces $U^\perp$ and $U^\perp \intersect p^\perp$ are transversal to all $F^j$. As a consequence we find that $\gen{U^\perp,p}$ is nearly transversal to all $F^j$, and that $U \intersect p^\perp$ is almost transversal to all $(F^j,\omega^j)$. Similarly, we see that $\gen{U^\perp \intersect p^\perp, p}$ is nearly transversal to all $F^j$, and that $\gen{U \intersect p^\perp,p} = \gen{U,p} \intersect p^\perp$ is almost transversal to all $(F^j,\omega^j)$.

So, by applying Lemma~\ref{lem:generalized_projection} to the pairs
$(U_i,U'_i)_{1 \le i \le 4} = ((U,U^\perp), (\gen{U,p},\gen{U^\perp,p}),$ $(U \intersect p^\perp,U^\perp \intersect p^\perp), (\gen{U \intersect p^\perp,p},\gen{U^\perp \intersect p^\perp,p}))$,
we obtain flags and forms $({F_i^j}',{\omega_i^j}')$ on $U'_i/U_i$ such that $W \pitchfork (F^j,\omega^j)$ if and only if $W \pitchfork ({F_i'}^j,{\omega_i'}^j)$ for $U_i < W < U'_i$.

Using the isomorphisms from Lemma \ref{lem:isomorphisms} we obtain flags $(F_i^j,\omega_i^j)$ on $U^\perp/U$ such that
\begin{enumerate}
\item $W/U \pitchfork (F_1^j,\omega_1^j)$ if and only if $W \pitchfork (F^j,\omega^j)$,
\item $W/U \pitchfork (F_2^j,\omega_2^j)$ if and only if $\gen{W,p}/\gen{U,p} \pitchfork ({F_2'}^j,{\omega_2'}^j)$ if and only if $\gen{W,p} \pitchfork (F^j,\omega^j)$,
\item $W/U \pitchfork (F_3^j,\omega_3^j)$ if and only if $(W \intersect p^\perp)/(U \intersect p^\perp) \pitchfork ({F_3'}^j,{\omega_3'}^j)$ if and only if $W \intersect p^\perp \pitchfork (F^j,\omega^j)$, and
\item $W/U \pitchfork (F_4^j,\omega_4^j)$ if and only if $\gen{W \intersect p^\perp,p}/\gen{U \intersect p^\perp,p} \pitchfork ({F_4'}^j,{\omega_4'}^j)$ if and only if $\gen{W \intersect p^\perp,p} \pitchfork (F^j,\omega^j)$.
\end{enumerate}

It remains to see that the flags are matching. If $j$ and $l$ are such that ${F^j}^\perp = F^l$, then ${F_1^j}^\perp = F_1^l$, ${F_4^j}^\perp = F_4^l$, ${F_2^j}^\perp = F_3^l$ and ${F_3^j}^\perp = F_2^l$.
\end{proof}

\section{Flip-flop systems and generalized Phan geometries}
\label{sec:flip-flop-gpg}

The aim of this section is to show that generalized Phan geometries arise quite naturally from flip-flop systems. Although it is possible to develop this theory over arbitrary fields, in view of our applications we will concentrate on flip-flop systems over finite fields of square order. This restriction allows us to take some shortcuts that would otherwise not be possible.

The setting is as follows. Let $q$ be a prime power, let $\Delta = (\Delta_+,\Delta_-,\delta_*)$ be a Moufang twin building regarded as a vertex-colored simplicial complex. Let $(W,S)$ be the type of $\Delta$, i.e., the associated Coxeter system. We assume that $\Delta$ arises from a group with an $\field_{q^2}$-locally split root group datum (see \cite{caprem09} and \cite[Chapter~8]{abrbro}) and say that $\Delta$ is \emph{defined over $\field_{q^2}$}.

\begin{dfn}
\label{dfn:flip}
A \emph{flip} $\theta$ is an involutory automorphism of $\Delta_+ \union \Delta_-$ such that
\begin{enumerate}
\item $\Delta_+^\theta = \Delta_-$ (and thus also $\Delta_-^\theta = \Delta_+$),
\item $\delta_\varepsilon(c,d)=\delta_{\bar{\varepsilon}}(c^\theta,d^\theta)$ for $(\varepsilon,\bar{\varepsilon}) \in \{(+,-),(-,+),(*,*)\}$, and\label{item:flips_are_isometries}
\item for every panel $\st p$ of $\Delta_+$ there exists a chamber $c \in \st p$ such that $\coproj_{\st p} c^\theta \ne c$.\label{item:flip_strong}
\end{enumerate}
\end{dfn}

Note that if $(W,S)$ is of spherical type, then the twin building $\Delta$ is just the twin building of the spherical building $\Delta_+$, cf.\ \cite[Example 5.136(a)]{abrbro}. In that case we also simply speak of a flip of $\Delta_+$.

The relatively new notion of a flip of a (twin) building is not uniform in the literature. We refer to \cite{horn08} for an overview. Flips as defined above were introduced in \cite{devmue07}.

For a chamber $c \in \Delta_+$ we denote by $\delta_\theta(c)$, the \emph{$\theta$-codistance of $c$}, the element $\delta_*(c,c^\theta)$ of $W$. Note that, if $R \le \Delta_+$ is a spherical residue parallel to its image $R^\theta$, then $(R,R^\theta,\delta'_*)$ is a twin building in a natural way. In this case $\delta'_*(c,c^\theta)$ equals $\delta_*(c,c^\theta)$ for $c \in R$ if and only if $R$ and $R^\theta$ are opposite in $\Delta$. 
A \emph{Phan chamber} with respect to a flip $\theta$ is a chamber $c$ with $\delta_\theta(c)=1$. The \emph{numerical $\theta$-codistance of $c$} is defined as $l_\theta(c) := l(\delta_\theta(c))$.

\begin{lem}
\label{lem:flips_are_strong}
Let $\theta$ be a flip. Let $c$ be a chamber and $w = \delta_\theta(c)$. If $s \in S$ is such that $l(sw) < l(w)$, then there exists a chamber $d$ that is $s$-adjacent to $c$ with $l_\theta(d) < l_\theta(c)$. In particular, the flip $\theta$ admits Phan chambers.
\end{lem}

\begin{proof}
Consider the $s$-panel $\st p$ of $c$. Note that $w$, being a $\theta$-codistance, is an involution, so that $sw = ws$. A chamber $d \ge p$ can have $\theta$-codistance $w$, $sw$ or $sws$ where $w$ is attained and either $l(sws) = l(w)-2$ or $sws=w$. So $\delta_*(\st p,\st p^\theta) = w$. By \eqref{item:flip_strong} of Definition \ref{dfn:flip} there is a chamber $d$ with $\coproj_{\st p} d^\theta \ne d$, so that in particular $\delta_\theta(d) \ne w$. Hence $l_\theta(d) < l_\theta(c)$. The last statement follows by induction.
\end{proof}

The following fact essentially states that every flip of a building of type $A_n$, $B_n$, or $C_n$ defined over a field is represented by a hermitian form.

\begin{prop}
\label{prop:phan_iff_non-degenerate}
Let $\field$ be a finite field. Let $V$ be a vector space over $\field$, equipped with a building geometry of type $A_n$, $B_n$ or $C_n$ and let $\theta$ be a flip. There is an involutory $\field$-automorphism $\sigma$ and a $\sigma$-hermitian form $\omega$ on $V$ such that $U$ is opposite $U^\theta$ if and only if $U$ is $\omega$-non-degenerate.

Moreover if $n = 1$, then we have the following: $\sigma = \id$ if the geometry contains at most $2$ non-Phan chambers and $\sigma$ is of order $2$ if it contains $q+1$ non-Phan chambers and $\field = \field_{q^2}$. It contains precisely $1$ non-Phan chamber only if $\chr \field = 2$.
\end{prop}

The proof consists of the following three lemmas. The first deals with the case where $n \ge 2$. The second and third cover the cases $A_1=C_1$ and $B_1$. In rank $1$ we consider the geometries as Moufang sets and flips are understood to be morphisms of Moufang sets (see \cite{demedsag09} and \cite{knop05}). As Moufang sets the types $A_1$, $B_1$, and $C_1$ all are isomorphic but in the concrete realization of a geometry of type $B_1$ producing a form is slightly more work.

\begin{lem}
Let $\field$ be a field, let $V$ be an $\field$-vector space. Let $n \ge 2$ and consider the geometry
\begin{enumerate}
\item of type $A_n$ of of all proper, non-trivial subspaces of $V$ where $\dim V = n+1$,
\item of type $B_n$ of all non-trivial subspaces of $V$ that are totally isotropic with respect to a non-degenerate symmetric bilinear form $(\cdot,\cdot)$ of Witt index $n$, where $\dim V = 2n+1$, or
\item of type $C_n$ of all non-trivial subspaces of $V$ that are totally isotropic with respect to a non-degenerate alternating bilinear form $(\cdot,\cdot)$, where $\dim V = 2n$.
\end{enumerate}
Let $\theta$ be a flip of the respective spherical (twin) building. There is an involutory $\field$-automorphism $\sigma$ and a $\sigma$-hermitian form $\omega$ on $V$ such that $U$ is opposite $U^\theta$ if and only if $U$ is $\omega$-non-degenerate.
\end{lem}

\begin{proof}
Let $w_0$ denote the longest word in $W$. That $\theta$ is a flip means that $\delta(c,d)=w_0\delta(c^\theta,d^\theta)w_0$ for two chambers $c,d$. Let $\tau$ denote the induced diagram automorphism. In the cases $B_n$ and $C_n$ the longest word $w_0$ is central so that $\delta(c,d) = \delta(c^\theta,d^\theta)$, hence $\tau = \id$. In the case $A_n$, $\tau$ is the non-trivial diagram automorphism. The image of the $J$-residue of $c$ under $\theta$ is the $\tau(J)$-residue of $\theta(c)$. In other words the induced map on the building geometry (that we will also denote by $\theta$) is type-preserving in the cases $B_n$ and $C_n$ and a polarity in the case $A_n$.

Case $A_n$: The polarity $\theta$ induces an isomorphism $\chi_\theta$ of the building geometry of $V$ with the building geometry of $V^*$ via $\chi(U) = (U^\theta)^\circ \defeq \{\alpha \in V^* \mid \alpha(U^\theta)=0\}$. By the Fundamental Theorem of projective geometry (see for example \cite[Theorem 2.26]{artin57}), there is an $\field$-automorphism $\sigma$ and a $\sigma$-semilinear map $\rho \colon V \to V^*$ that represents $\chi$. Now $U$ is opposite $U^\theta$ if and only if $U \intersect U^\theta = 0$ if and only if $U \intersect U^{\perp_\omega} = 0$ for the form $\omega'$ defined by $\omega'(x,y) = y^\rho(x)$. Since $\theta$ is involutory, $\omega'$ is reflexive. So by \cite[Theorem 6.1.3]{cameron91}, it is either alternating or a scalar multiple of a $\sigma$-hermitian form. If it were alternating then $U$ would never be opposite $U^\theta$. But we have seen in Lemma \ref{lem:flips_are_strong} that there is Phan chamber. So it must be a scalar multiple of a $\sigma$-hermitian form. Since we are just interested in whether or not a space is non-degenerate, scalar multiples do not play any role and we may take $\omega$ to be that $\sigma$-hermitian form.

Cases $B_n$, $C_n$: By \cite[Theorem 8.6, (II)]{tits74}, the automorphism $\theta$ can be extended to an automorphism of the projective geometry of $V$. This automorphism is again by the Fundamental Theorem of projective geometry represented by a $\sigma$-semilinear map $\rho$. Now $U$ is opposite $U^\theta$ if and only if $U \intersect (U^\theta)^\perp = 0$. So $U$ is opposite $U^\theta$ if and only if $U$ is $\omega$-non-degenerate where $\omega(x,y) = (x,y^\rho)$. As above one sees that $\omega$ may be chosen to be $\sigma$-hermitian.
\end{proof}

\begin{lem}
\label{lem:form_for_A_1}
Let $\field$ be a field and let $V$ be a two-dimensional $\field$-vector space. Let $\theta$ be an involutory automorphism of the Moufang set of the geometry of type $A_1 = C_1$ of $\prs(V)$. Assume that there is an $a \in V$ with $\gen{a}^\theta \ne \gen{a}$. There is an involutory $\field$-automorphism $\sigma$ and a $\sigma$-hermitian form $\omega$ such that a point is $\theta$-fixed if and only if it is $\omega$-isotropic.

Moreover there are $q+1$ points that are $\theta$-fixed if $\sigma \ne \id$ and $\field = \field_{q^2}$, there are two or no $\theta$-fixed points if $\theta = \id$ and $\chr \field \ne 2$ and there is one $\theta$-fixed point if $\theta = \id$ and $\chr \field = 2$.
\end{lem}

\begin{proof}
Let $b \in V$ be such that $\gen{b} = \gen{a}^\theta$. Let $X_1$ denote the Moufang set of $\prs(V)$ with $\infty = \gen{a}$ and $0 = \gen{b}$ and $X_2$ denote the Moufang set of $\prs(V)$ with $\infty = \gen{b}$ and $0 = \gen{a}$. Then $\theta$ can be seen as an isomorphism $X_1 \to X_2$ that fixes $0$ and $\infty$. So by \cite[Theorem 1]{hua49}, see also \cite[Theorem 3.3.1]{knop05}, there is an $\field$-automorphism $\sigma$ and scalars $\mu, \xi \in \field$ such that
\begin{equation}
\label{eq:theta_action}
\gen{\alpha a + \beta b}^\theta = \gen{\beta^\sigma \mu a + \alpha^\sigma \xi b} \textrm{ .}
\end{equation}
Consider the $\sigma$-hermitian form $\omega$ with $\omega(a,a) = -\mu$, $\omega(b,b) = \xi$ and $\omega(a,b) = 0$. We have $\omega(\alpha a + \beta b) = -\mu \alpha \alpha^\sigma + \xi \beta \beta^\sigma$ which vanishes if and only if
\begin{equation}
\label{eq:theta_fix}
\mu \alpha \alpha^\sigma = \xi \beta \beta^\sigma \textrm{ .}
\end{equation}
Comparing this with \eqref{eq:theta_action} (and recalling that $\theta$ interchanges $\gen{a}$ and $\gen{b}$) one sees that it is precisely the condition for $\gen{\alpha a + \beta b}$ to be $\theta$-stable.

The number of $\theta$-fixed points can be obtained by counting the number of solutions of \eqref{eq:theta_fix} (whether there are two or no solutions in the case $\chr \field \ne 2$, $\sigma = \id$ depends on whether $\mu/\xi$ is a square).
\end{proof}

\begin{lem}
\label{lem:form_for_B_1}
Let $\field$ be a field, let $V$ be a $3$-dimensional $\field$-vector space and let $(\cdot,\cdot)$ be a non-degenerate symmetric bilinear form on $V$ of Witt index $1$. Let $\theta$ define a flip of the Moufang set $X$ of the geometry of type $B_1$ defined by $V$ and $(\cdot,\cdot)$. Assume that there is a $g \in V$ such that $\gen{g}$ is a point of the geometry that is not $\theta$-fixed (i.e. $(g,g)=0$ and $\gen{g}^\theta \ne \gen{g}$). Then there is an involutory $\field$-isomorphism $\sigma$ and a $\sigma$-hermitian form $\omega$ such that a $(\cdot,\cdot)$-isotropic point is $\theta$-fix if and only if it is $\omega$-isotropic.

If $\field$ is finite, the correspondence between $\theta$-fixed points and the order of $\sigma$ is as in Lemma~\ref{lem:form_for_A_1}.
\end{lem}

\begin{proof}
Assume first that there is no $\theta$-fixed point. Let $y$ be a point that is $(\cdot,\cdot)$-anisotropic and let $H$ be a complement for $\gen{y}$ in $V$. Since $X$ is isomorphic to the Moufang set considered in Lemma \ref{lem:form_for_A_1}, we know that there is a $\sigma$-hermitian form on a two-dimensional $\field$-vector space that has no isotropic points. Let $\omega'$ be such a form on $H$. Define $\omega$ by letting $\omega(u,v) = \omega'(\pr_H(u),\pr_H(v))$ where $\pr_H$ is the projection onto $H$ that corresponds to the decomposition $V = H \oplus \gen{y}$. Now $\pr_H(v) \ne 0$ for every $(\cdot,\cdot)$-isotropic $v$, so $\omega$ is as desired.

So assume now that there is a $\theta$-fixed point $\gen{f}$, i.e. $(f,f) = 0$ and $\gen{f}^\theta = \gen{f}$. Let $e$ be $(\cdot,\cdot)$-isotropic such that $(e,f) = 1$. Let $x$ be $(\cdot,\cdot)$-orthogonal to $\gen{e,f}$. By normalizing $x$ and considering the form $\xi(\cdot,\cdot)$ if $\xi=(x,x)$ is a non-square we may assume that $\xi = 1$. Hence $(\cdot,\cdot)$ has Gram matrix
\[
\left(
\begin{array}{ccc}
&&1\\
&1&\\
1&&
\end{array}
\right)
\]
with respect to the basis $(e,x,f)$.

Let $H = \gen{e,x}$. The map $\rho \colon X \to \prs(H)$ that takes $\gen{f}$ to $\gen{x}$ and $\gen{u}$ to $\gen{u,f} \intersect H$ for $u$ isotropic, $u \ne f$, is an isomorphism of Moufang sets (cf. \cite[Theorem 3.3.5, Section 4.4]{knop05}, \cite[Section 6]{cameroncg}):

Indeed the map
\begin{align*}
\SL_2(\field) & \to  \SO_3(\field)\\
\left(
\begin{array}{cc}
a & b\\
c & d\\
\end{array}
\right)
& \mapsto 
\left(
\begin{array}{ccc}
a^2 & -2 a b & -2 b^2\\
-a c & a d+c b & 2 b d\\
-1/2 c^2 & c d & d^2
\end{array}
\right)
\end{align*}
is an homomorphism where $\SO_3$ is the group of elements of $\SL_3(\field)$ that preserve $(\cdot,\cdot)$ in the basis $(e,x,f)$ and the elements of $\SL_2(\field)$ are taken with respect to the basis $(e,x)$. This homomorphism induces isomorphisms of the root groups as desired.

Now $\theta$ induces an involutory automorphism of the Moufang set of $\prs(H)$ that does not fix $\rho(\gen{g})$. So by Lemma \ref{lem:form_for_A_1} there is an $\field$-automorphism $\sigma$ and a $\sigma$-hermitian form $\omega'$ on $H$ such that $\gen{u}$ is $\theta$-fix if and only if $\rho(\gen{u})$ is $\omega'$-isotropic.

We define the $\sigma$-hermitian form $\omega$ on $V$ by $\omega(u,v) = \omega'(\pr_H(u),\pr_H(v))$ where $\pr_H$ is the projection that corresponds to the decomposition $V = H \oplus \gen{f}$. In other words if $A$ is the Gram matrix for $\omega'$ with basis $(e,x)$, then
\[
\left(
\begin{array}{c|c}
A& \\
\hline
&0\\
\end{array}
\right)
\]
is the Gram matrix for $\omega$ with respect to $(e,x,f)$.

If $u \in H \setminus \gen{x}$, then $\rho^{-1}(\gen{u}) = \gen{u + \alpha f}$ for some $\alpha \in \field$. But
\[
\omega(u + \alpha f,u + \alpha f) = \omega(u,u)
\]
so $\gen{u + \alpha f}$ is isotropic if and only if $\gen{u}$ is isotropic if and only if $\gen{u + \alpha f}$ is $\theta$-fix. Moreover $\gen{f}$ is isotropic and $\theta$-fix. Hence $\omega$ is as desired.
\end{proof}

We now introduce flip-flop systems and show that they are isomorphic to generalized Phan geometries.

\begin{dfn}
Given a flip of $\Delta = (\Delta_+,\Delta_-,\delta_*)$ and a residue $R$ of $\Delta_+$ we denote by $R_\theta$ the complex of chambers $c$ of $R$ for which $l_\theta(c)$ is minimal. This is called the \emph{flip-flop system} of $R$.
\end{dfn}

If $R$ and $R^\theta$ are parallel, such that $(R,R^\theta,\delta_*')$ is a 
twin building with a flip $\theta|_{R \union R^\theta}$, then $R_\theta$ 
consists of the Phan chambers with respect to $\delta_*'$. If $R$ and 
$R^\theta$ are not parallel, the situation is described by the following 
proposition. The precise description of the relative links in Proposition 
\ref{prop:R_theta_characterization} below is due to Bernhard M\"uhlherr.

\begin{prop}
\label{prop:R_theta_characterization}
Let $R$ be a spherical residue of $\Delta_+$ and let $Q = \coproj_{R} R^\theta$. Then the map $\theta' = \coproj_R \circ \theta$ is a flip on the spherical twin building associated to $Q$.

Moreover a chamber $c$ of $R$ is in $R_\theta$ if and only if the following conditions are satisfied:
\begin{enumerate}
\item $c$ is opposite $Q$ in $R$, and
\item $x = proj_Q(c)$ is opposite $x^{\theta'}$ in $Q$.
\end{enumerate}
\end{prop}

\begin{proof}
First let us see that $\theta'$ is a flip on $Q$. Since $\theta$ is an isometry, we have $(\coproj_R c^\theta)^\theta = \coproj_{R^\theta} c$, and $\coproj_R \coproj_{R^\theta} c = \proj_Q c$. So $\theta'$ is involutory on $Q$.

Since $Q$ is spherical, in order to verify property \eqref{item:flips_are_isometries} it suffices to show that $\delta^Q_+(c,d)=\delta^Q_-(c^{\theta'},d^{\theta'})$ for $c,d \in Q$, i.e. that $\delta_+(c,d) = w_0 \delta_+(c^{\theta'},d^{\theta'}) w_0$ where $w_0$ is the longest word in the Coxeter group of $Q$. Now
\begin{align*}
\delta_+(c,d) & =  \delta_*(c,\coproj_{R^\theta}c) \delta_*(\coproj_{R^\theta}c,d)\\
 & =  \delta_*(c,\coproj_{R^\theta}c) \delta_-(\coproj_{R^\theta}c,\coproj_{R^\theta}d) \delta_*(\coproj_{R^\theta}d,d)\textrm{ .}
\end{align*}

Here $\delta_*(c,\coproj_{R^\theta} c)$ is the longest word of the Coxeter group of $Q$, because for any $e \in Q$ we have
\[
l_*(c,e) = l_*(c, \coproj_Q c) - l_-(\coproj_Q c,e) \text{ .}
\]
So for $l_*(c, \coproj_Q c)$ to be maximal, $l_-(\coproj_Q c,e)$ has to be maximal. Hence $\delta_*(c,\coproj_{R^\theta}c)$ and $\delta_*(\coproj_{R^\theta}d,d)$ both are the longest word in the Coxeter group of $Q^\theta$ which is the same as that of $Q$.

The last condition for $\theta'$ to be a flip is, that for every panel $\st p$ of $Q$ there is a chamber $c \in \st p$ such that $\proj_{\st p} c^{\theta'} \ne c$. But this follows immediately from the fact that $\proj_{\st p} \coproj_{R} = \coproj_{\st p}$.

\begin{figure}[!htb]
\mycenter{\includegraphics{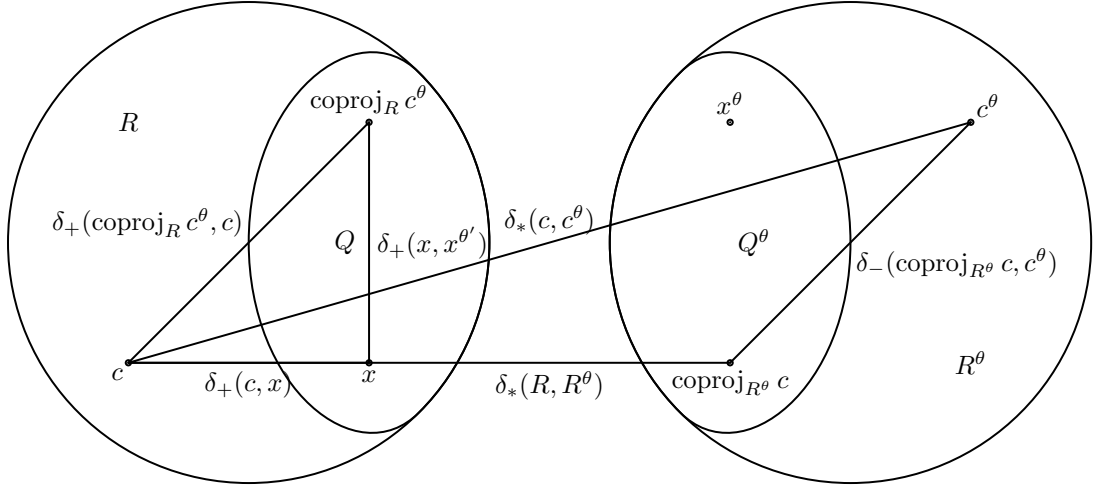}}
\caption{$l_*(c,c^\theta) = l_*(R,R^\theta) - 2l_+(c,x)-l_+(x,x^{\theta'})$}
\label{fig:R_theta_char}
\end{figure}

Now we want to prove the characterization of $R_\theta$. Let $c \in R$ be arbitrary. Then
\begin{align*}
\delta_\theta(c) & = \delta_*(c,c^\theta)\\
& = \delta_*(c,\coproj_{R^\theta}c) \delta_-(\coproj_{R^\theta}c,c^\theta)\\
& = \delta_+(c,\coproj_R \coproj_{R^\theta} c) \delta_*(\coproj_R \coproj_{R^\theta} c,\coproj_{R^\theta} c) \delta_-(\coproj_{R^\theta}c,c^\theta)
\end{align*}
(see Figure \ref{fig:R_theta_char}).
Here $\delta_*(\coproj_R \coproj_{R^\theta} c,\coproj_{R^\theta} c)$ is independent of the chamber $c$ and we denote it just by $\delta_*(R,R^\theta)$. Note that $\coproj_R \coproj_{R^\theta} c = \proj_Q c \eqdef x$. We have $\delta_-(\coproj_{R^\theta}c, c^\theta) = \delta_+(\coproj_R c^\theta,c)$ where $\coproj_R c^\theta = \coproj_Q c^\theta = x^{\theta'}$. And we can write $\delta_+(\coproj_R c^\theta,c)$ as $\delta_+(c,x)\delta_+(x,x^{\theta'})$. Putting this together we get that
\[
l_*(c,c^\theta) = l_*(R,R^\theta) - 2l_+(c,x)-l_+(x,x^{\theta'}) \textrm{ .}
\]
We see that this number is minimal if and only if $l_+(c,x)$ and $l_+(x,x^{\theta'})$ are both maximal among chambers $x \in Q$ and $c \in R$. For one can first choose $x$ so as to maximize $l_+(x,x^{\theta'})$ and then $c$ opposite $x$ inside $R$. And the longest word of $Q$ can be obtained as $l_+(x,x^{\theta'})$ by Lemma \ref{lem:flips_are_strong}.
\end{proof}

\begin{lem}
\label{lem:field-automorphisms_all_the_same}
Let $\Delta$ be a twin building of irreducible type defined over a finite field $\field_{q^2}$ and let $\theta$ be a flip of $\Delta$. The involutory $\field_{q^2}$-automorphism $\sigma$ whose existence is guaranteed by Proposition \ref{prop:phan_iff_non-degenerate} is the same for all spherical residues $R$ of types $A_n$, $B_n$, or $C_n$ that are parallel to their image under $\theta$.
\end{lem}

\begin{proof}
Since the field is finite, it suffices to distinguish between the cases $\sigma = \id$ and $\sigma \neq \id$. 
By the Lemmas~\ref{lem:form_for_A_1} and \ref{lem:form_for_B_1} it furthermore suffices to consider panels, because if $R$ is a residue of higher rank, we may take any Phan chamber $c$ of $(R,R^\theta,\delta_*')$, take a panel $P$ of $c$ and consider $P$ and $P^\theta$ instead of $R$ and $R^\theta$. It is clear, that the $\field_{q^2}$-automorphism for $\theta$ on $P$ and on $R$ is the same. 

So let $c$ be a chamber, let $s \in S$, and let $R$ be the $s$-panel of $c$. Let $(G,(U_\alpha)_{\alpha \in \Phi}, T)$ be the root group datum associated to $\Delta$ (see \cite[Section 8.5]{abrbro}). The flip induces an automorphism of the Kac--Moody group $G$ which we also denote by $\theta$ (see \cite{horn08}). By \cite[Corollary A]{capmue06}, since $\Delta$ is irreducible, $\theta$ splits as a product of an inner, a diagonal, a graph, a field and a sign automorphism. In particular, it acts with the same field automorphism on all subgroups it stabilizes. Let $L$ be a Levi factor of the stabilizer of the $s$-residue $R$ of $c$ and of $R^\theta$ (cf.\ \cite[Proposition 3.6]{capmue06}). Then $\theta|_L$ is an automorphism of $L$, since $\theta$ stabilizes $R \union R^\theta$. Note that a chamber $d \in R$ is in $R_\theta$ if and only if the root containing $d$ is not $\theta$-fixed.
The action of $L$ on $R$ is equivalent to the action of $\SL_2(\field)$ on a projective line. In \cite[Section 4]{demedgrahor07} the automorphisms of $\SL_2(\field)$ are described. It is also described that $\SL_2(\field)_\theta$ preserves a $\sigma'$-hermitian form where $\sigma'$ is the $\field$-automorphism of $\theta$. From this it is clear that the $\sigma'$ is of order $2$ if and only if the $\sigma$ in Lemma \ref{lem:form_for_A_1} is of order two.
\end{proof}

\begin{thm}
\label{thm:flip_flops_are_gpgs}
Let $\Delta$ be a twin building of irreducible type defined over a finite field $\field_{q^2}$ and let $\theta$ be a flip of $\Delta$. If $R$ is a spherical residue of $\Delta_+$ of type $A_n$, $B_n$, or $C_n$, then $R_\theta$ is isomorphic to a generalized Phan geometry of the respective type defined over $\field_{q^2}$.
\end{thm}

This theorem corresponds to \cite[Fact 5.1]{devgramue08} and extends it to residues of types $B_n$ and $C_n$.

\begin{proof}
The spherical residue $R$ is isomorphic to 
\begin{enumerate}
\item the geometry of proper, non-trivial subspaces of a vector space $V$ of dimension $N=n+1$ over $\field_{q^2}$, if $R$ is of type $A_n$;
\item the geometry of non-trivial, totally isotropic subspaces of vector space $V$ of dimension $N=2n+1$ over $\field_{q^2}$ that is equipped with a non-degenerate bilinear form of Witt index $n$, if $R$ is of type $B_n$; or
\item the geometry of non-trivial, totally isotropic subspaces of vector space $V$ of dimension $N=2n$ over $\field_{q^2}$ that is equipped with a non-degenerate alternating form, if $R$ is of type $C_n$.
\end{enumerate}

The residue $Q = \coproj_R R^\theta$ consists of chambers that contain a certain flag $F = (0=V_0 < \ldots < V_l=V)$. If $R$ is of type $B_n$ or $C_n$, we augment this flag to $(0 = V_0 < \ldots < V_l \le V_{l+1} \defeq V_l^\perp < \ldots < V_{2l+1} = V_0^\perp =V)$.

A chamber $C = (U_0 < \ldots < U_n)$ of $R$ is opposite $Q$ if and only if every $U_j$ is transversal to $F$ in case $A_n$, respectively $U_j^\perp$ is transversal to $F$ in the cases $B_n$ and $C_n$. To unify the cases let $\hat{C} = (U_0 < \ldots < U_n \le U_{n+1} \defeq U_n^\perp < \ldots < U_{2n+1} \defeq U_0^\perp)$ be the augmented chamber in $V$. The projection of $C$ to $Q$ is the chamber consisting of the spaces $\tgen{U_i,V_{k_{U_i}-1}} \intersect V_{k_{U_i}}  = \tgen{U_i \intersect V_{k_{U_i}},U_i}$ (respectively the isotropic part of this in cases $B_n$ and $C_n$).

Let $W_0 < \ldots < W_n$ be a chamber of $Q$ and $W_0' < \ldots < W_n'$ its image under $\theta$. The two are opposite in $Q$ if and only if for every $W_j$ the the following holds: Let $i$ be such that $V_{i-1} \le W_j < V_{i}$. Then $W_j/V_{i-1}$ has to be transversal to $W_j'/V_{i-1}$ in $V_{i}/V_{i-1}$ if the type is $A_n$ or $i \le l$. In the cases $B_n$ and $C_n$ for $i=l+1$ the condition becomes that ${W_j'}^\perp/V_l$ be transversal to $W_j/V_l$.

By Proposition \ref{prop:phan_iff_non-degenerate} there are forms $\omega_{i}'$ on $V_{i}/V_{i-1}$ such that $U/V_{i-1}$ is opposite $U^\theta/V_{i-1}$ if and only if $U$ is $\omega_{i}'$-non-degenerate for $V_{i-1} \le U \le V_{i}$. In the cases $B_n$ and $C_n$ this is only true for $i \le l+1$. For larger $i$ we define the $\omega_i'$ to make the square of natural isomorphisms commute:
\[
\begin{array}{ccc}
V_{i}/V_{i-1} & \stackrel{\omega_i'}{\simeq} & (V_{i}/V_{i-1})^*\\
\rsimeq {(\cdot,\cdot)} & & \rsimeq {(\cdot,\cdot)}\\
(V_{i-1}^\perp/V_{i}^\perp)^* & \stackrel{\omega_{2l+2-i}'}{\simeq} & V_{i-1}^\perp/V_i^\perp\text{ .} 
\end{array}
\]
The $\omega_i'$ induce forms $\omega_{i}$ on $V_{i}$ with radical $V_{i-1}$. By Lemma~\ref{lem:field-automorphisms_all_the_same} the $\field$-au\-to\-mor\-phisms $\sigma_i$ such that $\omega_i$ is $\sigma_i$-hermitian are all the same.

Proposition \ref{prop:R_theta_characterization} gives two conditions that have to be satisfied for $C$ to be in $R_\theta$. The first is, that $C$ is opposite $Q$. This is the case if and only if $U_j \pitchfork V_{k_{U_j}}$ for every $j$. The second is that the projection of $C$ to $Q$ is mapped to an opposite chamber in $Q$ under $\theta'$. If we translate this using our discussion above and set $i = k_{U_j}$, we obtain that $U_j \intersect V_i$ has to be $\omega_{i}$-non-degenerate.
\end{proof}

\begin{cor}
\label{cor:flip_flops_are_gpgs}
Let $\Delta = (\Delta_+,\Delta_-,\delta_*)$ be the affine twin building associated to $G=\SL_{n+1}(\field_{q^2}[t,t^{-1}])$ or to $G=\Sp_{2n}(\field_{q^2}[t,t^{-1}])$. Let $\theta$ be a flip of $\Delta$. If $R$ is an irreducible spherical residue of $\Delta_+$, then $R_\theta$ is isomorphic to a generalized Phan geometry of type $A_m$ or $C_m$ defined over $\field_{q^2}$.
\end{cor}

\begin{proof}
The residues of $\Delta_+$ are spherical buildings defined over $\field_{q^2}$. Since all irreducible residues of $\tilde{A}_n$ and $\tilde{C}_n$ are of type $A_m$ or $C_m$, the statement follows from the theorem.
\end{proof}

\section{Applications}
\label{sec:applications}

\subsection{Finiteness length of the unitary form of a Kac--Moody group}
\label{subsec:finiteness-length}

Let $\sigma$ be the involutory automorphism of the ring $\field_{q^2}[t,t^{-1}]$ that acts as the non-trivial involution on $\field_{q^2}$ and exchanges $t$ and $t^{-1}$. Let $G=\Sp_{2n}(\field_{q^2}[t,t^{-1}])$, let $\theta$ be the automorphism of $G$ that maps $g$ to $(g^\sigma)^{-T}$, and let $K:=G_\theta := \{ g \in G \mid g^\theta = g \}$.

The aim of this section is to show how the topological finiteness length of $K$ can be established using the Main Theorem.

\begin{thm}
\label{thm:finiteness_length}
The group $K$ is of type $\F_{n-1}$ but not of type $\FP_n$ provided $4^{n-1}(q+1)<q^2$.
\end{thm}

\begin{rem}
The group $K$ is an arithmetic subgroup of $\Sp_{2n}(\field_{q^2}(t))$, cf. \cite[Remark 5.6]{devgramue08}, of local rank $n$. So Theorem \ref{thm:finiteness_length} is an instance of the Rank Theorem, cf. \cite{buxgrawit10b}.
\end{rem}

Our main tool for establishing the finiteness length of $K$ is Brown's Criterion, see \cite{brown87}. We give a version that is particularly well suited for our purpose, it can be found in \cite{abramenko96} as Lemma 14.

\begin{prop}[Brown's criterion]
\label{prop:browns_criterion}
Let $X$ be a CW complex. Let $\Gamma$ act on $X$ by homeomorphisms that permute the cells. Assume that there exists an $n \ge 1$ such that the following conditions are satisfied
\begin{enumerate}
\item $X$ is $n$-connected. \label{item:n-connected}
\item If $\sigma$ is a cell of dimension $d \le n$, then the stabilizer $\Gamma_\sigma$ is of type $F_{n-d}$.\label{item:finite_stabilizers}
\item $X = \Union_{j \in \nat} X_j$ with $\Gamma$-invariant subcomplexes $X_j$ of $X$ which are finite complexes modulo $\Gamma$.\label{item:cocompact}
\item $X_{j+1} = X_j \union \Union_{a \in A_j} S_{a,j}$ with contractible subcomplexes $S_{a,j} \subseteq X_{j+1}$ satisfying \label{item:filtration_union}
\begin{enumerate}
\item $S_{a,j} \intersect S_{b,j} \subseteq X_j$ for all $j$ and $a \ne b \in A_j$,\label{item:intersection_contained}
\item $S_{a,j} \intersect X_j$ is $(n-1)$-spherical for all $j$ and all $a \in A_j$,\label{item:intersection_connected}
\item There exist infinitely many $j$ such that $\tilde{H}_{n-1}(S_{a,j} \intersect X_j) \ne 0$ for at least one $a \in A_j$. \label{item:homology_nonzero}
\end{enumerate}
\end{enumerate}
Then $\Gamma$ is of type $\F_{n-1}$ and not of type $\FP_n$.
\end{prop}

Recall that the group $\Sp_{2n}(\field_{q^2}[t,t^{-1}])$ admits a twin BN-pair that gives rise to a twin building $(\Delta_+,\Delta_-,\delta_*)$ (see \cite[Section 6.12]{abrbro} for the $A_n$ case). The group $G$ acts on the twin building by isometries, that is, $\delta_\varepsilon(c,d)=\delta_{{\varepsilon}}(gc,gd)$ for $\varepsilon \in \{+,-,*\}$ and $g \in G$.

We use the following lemma, which is taken from \cite{horn08}.

\begin{lem}
The automorphism $\theta$ of $G$ induces a flip of $(\Delta_+,\Delta_-,\delta_*)$.
\end{lem}

In our case, where the residue field is finite and $\theta$ induces the non-trivial field involution on the residue field, the flip $\theta$ also satisfies the following ascending version of Lemma \ref{lem:flips_are_strong}. (Recall from the preceding section that the $\theta$-codistance $\delta_\theta$ of $c$ is $\delta_*(c,c^\theta) \in W$ and that the numerical $\theta$-codistance of $c$ is the integer $l(\delta_\theta(c))$.)

\begin{lem}
\label{lem:flip_is_costrong}
For every panel $\st p$ of $\Delta_+$ there exists a chamber $c \ge p$ such that $\coproj_{\st p} c^\theta = c$. Thus if $c$ is a chamber, $\st p$ is its $s$-panel, $w = \delta_\theta(c)$, and $s \in S$ is such that $l(sw) > l(w)$, then there is a chamber $d \ge p$ such that $l_\theta(d) > l_\theta(c)$.
\end{lem}

\begin{proof}
If $\coproj_{\st p} \st p^\theta$ consists of only one chamber, there is nothing to show -- this chamber is as desired. So we may assume that $\coproj_{\st p} \st p^\theta = \st p$. Then $\theta' \defeq \coproj_{\st p} \circ \theta$ is a flip on $\st p$. Lemmas~\ref{lem:form_for_A_1} and \ref{lem:form_for_B_1} show that there is a chamber that is $\theta'$-fixed, i.e. satisfies $\coproj_{\st p} c^\theta = c$.
The second assertion is shown as in the proof of Lemma \ref{lem:flips_are_strong}.
\end{proof}

Note that $(kc)^\theta = k^\theta c^\theta = k c^\theta$, so that $\delta_\theta(kc) = \delta_\theta(c)$ for $k \in K$. In particular $K$ acts on the subcomplexes $\Delta_j=\overline{\{c \in \calC(\Delta_+) \mid l_\theta(c) \le j\}}$ of $\Delta_+$. We claim that the filtration $(X_j)_{j \in \nat} = (\trealize{\Delta_j})_{j \in \nat}$ and the space $X = \realize{\Delta_+}$ satisfy the conditions of Brown's Criterion (Proposition \ref{prop:browns_criterion}).

We take $A_j$ to be the set of simplices which are in $\Delta_{j+1}$ but not in $\Delta_j$ and all of whose proper faces are in $\Delta_j$. Moreover, we let $S_{a,j} = \trealize{\st_{\Delta_{j+1}} a}$.

\begin{lem}
\label{lem:chamber_contains_unique_a_i}
We have $X_{j+1} = \Union_{i \in A_j} X_j$. More precisely every chamber $c \in \Delta_{j+1}$, $c \nin \Delta_j$ has a unique face that is in $A_j$.
\end{lem}

\begin{proof}
First let us see that, if $b \le c$ is a face that is contained in $\Delta_j$, then there is a facet (maximal proper face) $p$ of $c$ that contains $b$ such that $p$ is contained in $\Delta_j$. Let $d \ge b$ be a chamber with $l_\theta(d) < l_\theta(c)$. By \cite{horn08}, there is a chamber $c' \ge b$ that is adjacent to $c$ with $l_\theta(c') < l_\theta(c)$. So $p = c' \intersect c$ is as desired.

So $\overline{c} \intersect \Delta_j$ consists of (rather, is the complex generated by) facets $p_1, \ldots, p_k$ of $c$. Let $v_1,\ldots,v_k$ be the vertices of $c$ such that $v_i \nin p_i$ and let $a$ be the simplex spanned by $v_1, \ldots, v_k$. Clearly $a$ is in $A_j$. Conversely if $a' \in A_j$, $a' \le c$, then for every facet $p_j$ there is a vertex $v_j$ of $a'$ not contained in $p_j$ (because otherwise $a' \le p_j$). If there were another vertex, say $v$, then $a' \setminus v$ would not be in $\Delta_j$, hence $a = a'$.
\end{proof}

\begin{lem}
\label{lem:st_intersect_Xj}
Let $j \in \nat$, $a \in A_j$ and set $R = \lk_{\Delta_+} a$. Then $\lk_{\Delta_{j+1}} a = R_\theta$ and $\st_{\Delta_{j+1}} a \intersect \Delta_{j} = \lk_{\Delta_{j+1}} a * \partial a$.
\end{lem}

\begin{proof}
Note that chambers $c \in \st_{\Delta_+} a$ satisfy $l_\theta(c) \ge j+1$ and are in $\st_{\Delta_{j+1}} a$ if and only if $l_\theta(c) = j+1$. So these are precisely the chambers that have minimal numerical $\theta$-codistance among the chambers of $\st_{\Delta_{j+1}} a$. This shows the first statement.

For the second we have to see that each chamber $c \in \st_{\Delta_{j+1}} a$ is adjacent to a chamber of strictly shorter numerical $\theta$-codistance along every facet that does not contain $a$. Let $c$ be as above and let $\st p$ be a panel that does not contain $a$. Then $p \intersect a$ is a facet of $a$, so it is contained in $\Delta_j$. By the proof of the preceding lemma there is a panel of $c$ in $\Delta_j$, that contains $p \intersect a$. This panel, of course, cannot contain $a$, so it is $p$.
\end{proof}

\begin{lem}[cf. {\cite[Lemma 2.4]{abels91}}]
\label{lem:disjoint_or_identical}
Let $j \in \nat$ and $a,a' \in A_j$. Either $a = a'$ or $\st_{\Delta_{j+1}} a \intersect \st_{\Delta{j+1}} a' \subseteq \Delta_j$.
\end{lem}

\begin{proof}
By the preceding lemma if $b \in \st_{\Delta_{j+1}} a \intersect \st_{\Delta{j+1}} a'$, but $b \nin \Delta_j$, then $b \ge a$ and $b \ge a'$. Let $c$ be a chamber containing $b$. By Lemma \ref{lem:chamber_contains_unique_a_i} there is a unique face of $c$ that is one of the $a$, hence $a = a'$.
\end{proof}

\begin{lem}
\label{lem:properly_spherical}
For every $j \in \nat$ there is a $j' \ge j$, a chamber $c \in \Delta_+$ and an $s \in S$ such that $w = \delta_\theta(c)$ satisfies that $l(w) = j'$, and $l(tw)<l(w)$ for all $t \in S \setminus \{s\}$. Moreover there is a chamber $c'$ that is $s$-adjacent to $c$ with $\delta_\theta(c') = w$.
\end{lem}

\begin{proof}
Let $s \in S$ be arbitrary. We start with a chamber $c_0$ with $l_\theta(c_0) \ge j$, which is possible by Lemma \ref{lem:flip_is_costrong}. As long as possible, we take $t$ from $S \setminus \{s\}$ such that $l(tw_j) > l(w_j)$. By Lemma \ref{lem:flip_is_costrong}, there is a chamber $c_{j+1}$ with $\delta_\theta(c_{j+1}) = tw_j$ or $\delta_\theta(c_{j+1}) = tw_jt$ and $l_\theta(c_{j+1}) > l_\theta(c_j)$. We let $w_{j+1} = \delta_\theta(c_{j+1})$. This process has to terminate because $S \setminus \{s\}$ is spherical, so at some point $w \defeq w_j$ begins with the longest word of $S \setminus \{s\}$. The chamber $c \defeq c_j$ is as desired, see Lemma \ref{lem:flips_are_strong}.

Now clearly $l(s w) > l(w)$, because otherwise $w$ would be a longest word of $(W,S)$ which does not exist. We consider the $s$-panel $\st p$ of $c$. Either all chambers in $\st p$ have $\theta$-codistance $w$ or $sw$ or they all have $\theta$-codistance $w$ or $sws$. If $\coproj_{\st p}(\st p^\theta)$ consists of just one chamber, any other chamber will have smaller numerical $\theta$-codistance, i.e. $\theta$-codistance $w$. Since the building is thick, there is at least one which is distinct from $c$. We pick one and take it to be $c'$. If $\coproj_{\st p}(\st p^\theta) = \st p$, we may take $c'$ to be $\coproj_{\st p} c^\theta$.
\end{proof}

For the proof of Theorem \ref{thm:finiteness_length} we need the Lang--Steinberg Theorem, which is Theorem~J in \cite{steinberg77}:

\begin{thm}
\label{thm:lang}
Let $\bfG$ be a connected algebraic group defined over $\field_q$. Let $\theta$ be an endomorphism of $G \defeq \bfG(\overline{\field_q})$. If $G_\theta$, the set of $\theta$-fixed points of $G$, is finite, then the map $g \mapsto g^{-\theta}g$ is surjective.
\end{thm}

\begin{proof}[Proof of Theorem \ref{thm:finiteness_length}]
We verify the prerequisites of Brown's Criterion (Proposition \ref{prop:browns_criterion}). The first is clear by the Solomon-Tits Theorem (\cite[Theorem 1]{solomon69}, see also \cite[Section IV.6]{brownb}, \cite[Section 4.12]{abrbro}).

Now we want to verify condition \eqref{item:finite_stabilizers}. In fact we argue that the stabilizer $\stab_K(a)$ is finite for any (non-empty) simplex $a$ of $\Delta_+$, which implies that $\stab_K(a)$ is of finiteness type $\F_\infty$. 

First we consider the case that $a = c$ is a chamber. To do so, we use facts from \cite[Section 8.2]{abrbro}. Fix a twin apartment $\Sigma=(\Sigma_+,\Sigma_-)$ that contains $c$ and $c^\theta$.
By Proposition 8.15 and Proposition 8.19 of \cite{abrbro} the stabilizer of $c \union c^\theta$ is $U_{\alpha_1} \cdots U_{\alpha_m} H$ where $\alpha_1,\ldots,\alpha_m$ is a certain ordering of the twin roots $\alpha$ of $\Sigma$ with $c,c^\theta \in \alpha$, the group $U_\alpha$ is the root group of $\alpha$, and $H$ is the torus. So $\stab_G(c \union c^\theta)$ is finite, since $H$ and every of the $U_\alpha$ is. Hence $\stab_K(c)$ is finite.

Now let $a \in \Delta_+$ be arbitrary and let $c \ge a$ be a chamber. Once we realize that $\stab_{\stab_K(a)}(c)$ is finite and that $c$ has finite orbit under $\stab_K(a)$, finiteness of $\stab_K(a)$ follows from the orbit-stabilizer formula. That $\stab_{\stab_K(a)}(c)$ is finite is clear from the above argument because $\stab_{\stab_K(a)}(c) \le \stab_K(c)$. That the orbit of $c$ under $\stab_K(a)$ is finite follows from the fact that $\Delta_+$ is locally finite.

We claim that $K$ acts transitively on the chambers that have a given $\theta$-codistance. This immediately implies \eqref{item:cocompact}. Let $\bfT$ be the torus of diagonal matrices as an algebraic group scheme defined over the residue field $\field_{q^2}$. Let $(\Sigma_+,\Sigma_-)$ be the twin apartment corresponding to $\bfT$. Let $c'$ be an arbitrary chamber and let $c \in \Sigma_+$ be such that $c^\theta \in \Sigma_-$ and $\delta_\theta(c)=\delta_\theta(c')$. Let $(\Sigma_+',\Sigma_-')$ be a twin apartment such that $c' \in \Sigma_+'$,  ${c'}^\theta \in \Sigma_-'$. By strong transitivity there is a $g \in G$ such that $g \Sigma_+ = \Sigma_+'$ and $g \Sigma_- = \Sigma_-'$ and there is an $h \in G$ normalizing $(\Sigma_+',\Sigma_-')$ such that $hgc=c'$. Then $hgc^\theta$ is in $\Sigma_-'$ and has Weyl-distance $\delta_\theta(c')$ from $c'$, so $hgc^\theta={c'}^\theta$. Now
\[
(hg)^{-\theta}hgc = (hg)^{-\theta}c' = ((hg)^{-1}{c'}^\theta)^\theta = c
\]
and, similarly, $(hg)^{-\theta}hgc^\theta=c^\theta$, so $t = (hg)^{-\theta}hg \in T(\field_{q^2})$. Moreover $t^\theta = (hg)^{-1}(hg)^\theta = t^{-1}$.

Let $\overline{\field_{q^2}}$ denote the algebraic closure of $\field_{q^2}$ and note that $T(\overline{\field_{q^2}})$ is connected with respect to the Zariski topology. Let $\sigma$ be the endomorphism of raising elements of $\overline{\field_{q^2}}$ to the $q$th power and note that the fixed points set of $\overline{\field_{q^2}}$ under $\sigma^2$ is exactly $\field_{q^2}$. The map $\theta \colon \bfT(\overline{\field_{q^2}}) \to \bfT(\overline{\field_{q^2}}), g \mapsto g^{-\sigma}$ is an endomorphism, which satisfies $g^{\theta^2}=g^{\sigma^2}$, so that $s \in \bfT(\field_{q^2})$ if and only if $s=s^{\theta^2}$ for $s \in \bfT(\overline{\field_{q^2}})$.

Consequently $\bfT(\overline{\field_{q^2}})_{\theta} \subseteq \bfT(\overline{\field_{q^2}})_{\theta^2} = \bfT(\field_{q^2})$ is finite and so by the Lang--Steinberg Theorem (Theorem~\ref{thm:lang}) there is an $s \in \bfT(\overline{\field_{q^2}})$ such that $s^{-\theta}s=t$. Now
\[
s^{\theta^2} = (s t^{-1})^\theta =  s^\theta t = s \textrm{ ,}
\]
so $s \in \bfT(\field_{q^2})$. Hence the element $hgs$ lies in $K$ and it maps $c$ to $c'$. Therefore $K$ acts transitively on the chambers that have a given $\theta$-codistance.

Property \eqref{item:intersection_contained} follows from Lemma \ref{lem:disjoint_or_identical}.
As for \eqref{item:intersection_connected} we know by Lemma \ref{lem:st_intersect_Xj} that $S_{a,j} \intersect X_j = \trealize{(\lk_{\Delta_+} a)_\theta} * \trealize{\partial a}$, where $\lk_{\Delta_+} a$ is a residue of type $A_n$ or $C_n$. If $a$ has dimension $d$, then by Corollary \ref{cor:flip_flops_are_gpgs} $(\lk_{\Delta_+} a)_\theta$ is isomorphic to a generalized Phan geometry of type $A_{n-d}$ or $C_{n-d}$. Using the Main Theorem we get that $(\lk_{\Delta_+} a)_\theta$ is $(n-d-1)$-spherical in the $C_n$ case. In the $A_n$ case apply the Main Theorem of \cite{devgramue08}. Since $\partial a$ is a $(d-1)$-sphere, $\trealize{(\lk_{\Delta_+} a)_\theta} * \trealize{\partial a}$ is $(n-1)$-spherical by Proposition~\ref{prop:topology_basics}~\eqref{item:join_of_sphericals}.

The condition \eqref{item:homology_nonzero} follows from Lemma \ref{lem:properly_spherical} as follows: Let $j \in \nat$ be arbitrary and let $c$, $c'$, $j'$ and $s$ be as in the Lemma. Let $\st p$ be the $s$-panel of $c$. Then $p$ is in $A_{j'}$, because for every $t \in S \setminus \{s\}$ by Lemma \ref{lem:flips_are_strong} there is a chamber $d$ that is $t$-adjacent to $c$ with $l_\theta(d) < l_\theta(c)$. Now $(\lk_{\Delta_+} p)_\theta$ contains at least two vertices, namely $c \setminus p$ and $c' \setminus p$. So $\trealize{(\lk_{\Delta_+} p)_\theta} * \trealize{\partial p}$ is the join of an $(n-2)$-sphere with a space that is properly $0$-spherical, hence properly $(n-1)$-spherical.
\end{proof}

The proof of Theorem~\ref{thm:finiteness_length} also applies to Kac--Moody groups of hyperbolic type provided all proper residues are of type $A_m$ or $C_m$. To make this precise let $\bfG$ be a Kac--Moody group, and let $\sigma$ be the involutory $\field_{q^2}$-automorphism. There is a unique automorphism $\theta$ of $G \defeq \bfG(\field_{q^2})$ which on the defining rank-$1$-subgroups acts as $g \mapsto (g^\sigma)^{-T}$. Let $K \defeq G_\theta$ be the group of fixed points of this involution.

\begin{figure}
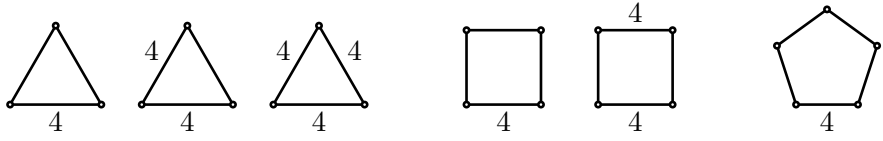

\begin{center}
\includegraphics{hyperbolic_diagrams-3}
\hspace{1cm}
\includegraphics{hyperbolic_diagrams-4}
\hspace{1cm}
\includegraphics{hyperbolic_diagrams-5}
\end{center}
\caption{Hyperbolic Coxeter diagrams whose proper subdiagrams are spherical.}
\label{fig:hyperbolic_diagrams}
\end{figure}

\begin{thm}
\label{thm:finiteness_length_hyperbolic}
Let $K$ be as above where the Coxeter diagram associated to $\bfG$ is one those shown in Figure~\ref{fig:hyperbolic_diagrams} and has $n$ vertices. Assume that $4^{n-1}(q+1)<q^2$.\\\noindent
If $n=3$, then $K$ is (finitely generated but) not of type $\FP_2$.\\\noindent
If $n=4$, then $K$ is finitely presented but not of type $\FP_3$.\\\noindent
If $n=5$, then $K$ is finitely presented.
\end{thm}

\begin{proof}
The first to cases are proved as Theorem~\ref{thm:finiteness_length}. The last case follows with \cite[Theorem~3.14]{devmue07} and \cite[Corollaire~1]{tits86} from Theorem~\ref{thm:geometry_simply_connected} below.
\end{proof}


\subsection{Local recognition of groups}
\label{subsec:local_recognition}

\begin{thm}
\label{thm:geometry_simply_connected}
Let $\field_{q^2}$ be a finite field of square order and let $\sigma$ be the non-trivial $\field_{q^2}$-involution. Let $X$ be an irreducible diagram of rank at least $3$, all irreducible rank $3$-residues of which are of type $A_3$, $B_3$, or $C_3$. Assume that $q^2 \geq 16(q+1)$. Let $\Delta$ be a twin building of type $X$ defined over $\field_{q^2}$ with a flip $\theta$ that induces $\sigma$ on $\field_{q^2}$ in the sense of Lemma \ref{lem:field-automorphisms_all_the_same}.

Then $\Delta_\theta=\overline{\{c \in \calC(\Delta_+) \mid \delta_\theta(c) = 1 \}}$ is $2$-simply connected.
\end{thm}

\begin{proof}
The geometry of $\Delta_+$ is $2$-simply connected, see \cite[Theorem 3]{tits81}. We want to use (the proof of) \cite[Theorem 3.14]{devmue07} to see that $\Delta_\theta$ is $2$-simply connected if and only if $\Delta_+$ is. To do so, we have to show that the residues of rank $3$ are simply connected and the residues of rank $2$ are connected. But this is true by Theorem \ref{thm:flip_flops_are_gpgs} and Theorem \ref{thm:precise_statement}, respectively \cite[Main Theorem]{devgramue08}.
\end{proof}

The following result has been published in \cite{bengrahofshp07}, \cite{grahornic07} in the case $B_n$ and \cite{grahofshp03}, \cite{grahornic06} in the case $C_n$. For the $F_4$ case Hoffman, Köhl, Mühlherr, and Shpectorov found a proof in Oberwolfach in summer 2005, which because of its length and its tedious case distinctions they did not publish. Our proof is much shorter as the concept of generalized Phan geometries allowed us to get rid of all case distinctions except the distinction of the Dynkin diagrams. 

\begin{thm}[Hoffman, Köhl, Mühlherr, Shpectorov 2005; unpublished]
\label{thm:local_recognition}
Let $n \ge 3$, let $q \ge 17$, let $X \in \{B_n,C_n,F_4\}$, and let $K$ be a group with a weak Phan system of type $X$ over $\field_{q^2}$. Then $K$ is a central quotient of $\Spin_{2n+1}(q)$, $\Sp_{2n}(q)$, or the simply connected version of $F_4(q)$, respectively.
\end{thm}

\begin{proof}
By Theorem \ref{thm:geometry_simply_connected}, $\Delta_\theta$ is $2$-simply connected, so in particular it is simply connected. Tits' Lemma \cite[Corollaire 1]{tits86} then implies that $K$ is the universal enveloping group of the weak Phan amalgam in $K$, cf. \cite{dunlap05}. The amalgams that can occur have been classified in \cite{benshp04}, \cite{gramlich04}, \cite{dunlap05} and they are unique up to passage to quotients.
\end{proof}

\begin{rem}
Theorem~\ref{thm:local_recognition} also holds for all unitary forms $K$ of arbitrary three-spherical tree type because spherical residues of rank three are either reducible or of type $A_3$, $B_3$, or $C_3$. If the three-spherical diagram is not a tree, the Phan system is not uniquely determined, but $K$ nevertheless is the universal enveloping group of the Phan system it contains.
\end{rem}

\bibliographystyle{amsalpha}
\bibliography{gpg_ams}

\end{document}